\documentclass[5p,12pt]{elsarticle}

\usepackage{bm}
\usepackage{amsmath}
\usepackage{amssymb}
\usepackage{mathtools}
\usepackage{siunitx}
\usepackage{fancyref}
\usepackage{pgfplots}
\usepackage{pgfplotstable}
\usepgfplotslibrary{external}

\usepackage{subcaption}
\usepackage{todonotes}
\usepackage{eurosym}
\usepackage{multirow}
\usepackage{tabularx}
\usepackage{booktabs}
\usepackage{subcaption}
\usepackage{mwe}
\usepackage[nomessages]{fp}
\usepackage{booktabs}
\usepackage{tabularx}

\usepackage{todonotes}

\usepackage[acronym]{glossaries}

\usepackage{layout}
\usepackage{printlen}
\usepackage{hyperref}

\usepackage{placeins} 
\usepackage{tablefootnote} 

\usepackage{nomencl} 
\makenomenclature 
\usepackage{etoolbox} 
\usepackage{url} 

\usetikzlibrary{math} 

\DeclareSIUnit{\ct}{ct}
\DeclareSIUnit\year{yr}

\newacronym{nlp}{NLP}{Non Linear Problem}
\newacronym{milp}{MILP}{Mixed Integer Linear Program}
\newacronym{minlp}{MINLP}{Mixed Integer Non-Linear Program}
\newacronym{dhn}{DHN}{District Heating Network}
\newacronym{dhns}{DHNs}{District Heating Networks}
\newacronym{gis}{GIS}{Geographic Information System}
\newacronym{chp}{CHP}{Combined Heat and Power}
\newacronym{4gdh}{4GDH}{4th Generation District Heating}
\newacronym{capex}{CAPEX}{Capital Expense}
\newacronym{opex}{OPEX}{Operational Expense}

\newcommand{\COTWO}{\mathrm{CO_2}}
\newcommand{\COTWOITALIC}{CO_2}
\newcommand{\FGDH}{\mathrm{4GDH}}
\newcommand{\FG}{\mathrm{4G}}

\DeclareSIUnit{\sieuro}{\mbox{\euro}}

\makeatletter
\renewcommand{\todo}[2][]{\tikzexternaldisable\@todo[#1]{#2}\tikzexternalenable}
\makeatother

\begin{document}
\emergencystretch 3em 

\setlength{\parindent}{0cm}


\begin{frontmatter}

\title{Decarbonization of Existing Heating Networks through Optimal Producer Retrofit and Low-Temperature Operation}

\author[1,3]{Martin Sollich} 
\ead{martin.sollich@kuleuven.be}
\author[1,2,3]{Yannick Wack}
\author[2,3]{Robbe Salenbien}
\author[1,3]{Maarten Blommaert}

\affiliation[1]{organization={Department of Mechanical Engineering, KU Leuven},addressline={ Celestijnenlaan 300 box 2421},city={Leuven},postcode={3001},country={Belgium}}

\affiliation[2]{organization={Flemish Institute for Technological Research (VITO)},addressline={Boeretang 200},city={Mol},postcode={2400},country={Belgium}}

\affiliation[3]{organization={EnergyVille},addressline={Thor Park, Poort Genk 8310},city={Genk},postcode={3600},country={Belgium}}

\begin{abstract}
District heating networks are considered a key factor for enabling emission-free heat supply, while many existing district heating networks still heavily rely on fossil fuels. With district heating network pipes easily exceeding a lifetime of 30 years, there is a growing potential to retrofit the heat producers of existing networks to enable low-emission heat supply. Today, the heat producer retrofit for district heating networks usually focuses on simplified approaches, where the non-linear nature of the design problem is relaxed or not considered at all. Some approaches take non-linearities into account but use optimization routines that are either not scalable to large problems or are not reliable in obtaining an optimal solution, such as parameter optimization and sensitivity studies. This paper presents an automated design approach, to decarbonize existing heating networks through optimal producer retrofit and ultimately enabling 4th generation operation. The approach uses multi-objective, mathematical optimization to balance $\COTWO$ emissions and network costs, by assessing different $\COTWO$ prices, and is based on a detailed physical model. The optimizer is given the freedom to choose the producer types, their capacities, and for each period, their supplied heat and supply temperature. The heat producers considered in this study are a natural gas boiler, an air-source heat pump, a solar thermal collector, and an electric boiler. A non-linear heat transport model accurately accounts for heat and momentum losses throughout the network, and ensures the feasibility of the proposed design and operation. The multi-period formulation incorporates temporal changes in heat demand and environmental conditions throughout the year. By formulating a continuous problem and using adjoint-based optimization, the automated approach remains scalable towards large scale applications. The design approach was assessed on a medium-sized 3rd DHN case and was able to optimally retrofit the heat producers. The retrofit study highlights a strong influence of the $\COTWO$ price on the optimal heat producer design and operation. Increasing $\COTWO$ prices shift the design toward a heat supply dominated by an energy-efficient and low-carbon heat pump. Furthermore, it was observed that even for the highest explored $\COTWO$ price of 0.3$\si{\,\sieuro}\,\mathrm{kg^{-1}}$, the low-carbon heat pump, electric boiler and solar thermal collector can not fully replace the natural gas boiler in an economic way. 
\end{abstract}


\end{frontmatter}

\setcounter{footnote}{0}

\newcommand{\tp}[1]{#1^{\intercal}} 			
\DeclarePairedDelimiter\abs{\lvert}{\rvert} 	
\newcommand{\card}[1]{\lvert#1\rvert}
\newcommand{\infNorm}[1]{\|#1\|_{\infty}}
\newcommand{\Real}[1]{\mathbb{R}^{#1}}
\newcommand{\ve}[1]{\bm{#1}} 	

\newcommand{\equalNew}{\tilde{\equalCon}}
\newcommand{\inEqualNew}{\tilde{\inEqualCon}}
\newcommand{\designVarNew}{\tilde{\designVar}}
\newcommand{\designVarTime}{\ve{\varphi}_{\timeVar}}

\newcommand{\radFlowNew}{\tilde{\radValve}}
\newcommand{\prodInputNew}{\tilde{\prodInput}}
\newcommand{\equalConModel}{\equalCon_\mathrm{m}}
\newcommand{\inEqualModel}{\equalCon_\mathrm{s}}
\newcommand{\stateVarModel}{\stateVar_\mathrm{m}}
\newcommand{\stateVarIneq}{\stateVar_\mathrm{s}}

\newcommand{\ALagrangian}{\mathcal{L}}
\newcommand{\LagMultis}{\lambda}
\newcommand{\LagPen}{\mu}
\newcommand{\slack}{s}
\newcommand{\equalConState}{g}

\newcommand{\costFull}{\mathcal{J}}
\newcommand{\costi}[1]{\cost_{\mathrm{#1}}}

\newcommand{\subCAPEX}{CAP}
\newcommand{\subOPEX}{OP}
\newcommand{\subPipe}{pipe}
\newcommand{\subHeat}{h}
\newcommand{\subCOTWO}{CO_2}
\newcommand{\subPump}{p}
\newcommand{\subRev}{rev}
\newcommand{\periodWeights}{w}
\newcommand{\Jpipepol}{\kappa}
\newcommand{\Jpipesmooth}{\xi}
\newcommand{\cPipe}{\npvCost_{mathrm{\subPipe}}}

\newcommand{\cHeatOPEX}{\npvCost_{\mathrm{hO}}}
\newcommand{\cHeatOPEXi}{\npvCost_{\mathrm{O,prod}}}
\newcommand{\cHeatOPEXiIJ}{\npvCost_{\mathrm{O,prod},ij}}
\newcommand{\cHeatOPEXiIJbold}{\bm{\npvCost_{\mathrm{O,prod},ij}}}

\newcommand{\KOPEX}{K}
\newcommand{\KOPEXvalue}{8208} 

\newcommand{\pumpEff}{\efficiency_{\mathrm{pump}}}
\newcommand{\cPumpOPEX}{\npvCost_{\mathrm{pO}}}
\newcommand{\cPumpOPEXi}{\npvCost_{\mathrm{O,pump}}}
\newcommand{\cPumpCAPEX}{\npvCost_{\mathrm{pC}}}
\newcommand{\cPumpCAPEXi}[1]{\npvCost_{\mathrm{pC},#1}}

\newcommand{\cRev}{\npvCost_{\mathrm{r}}}
\newcommand{\cRevi}[1]{\npvCost_{\mathrm{r},#1}}
\newcommand{\flow}{q}
\newcommand{\pressure}{p}
\newcommand{\temperature}{T}

\newcommand{\Reynolds}{Re}      
\newcommand{\density}{\rho}     
\newcommand{\viscosity}{\mu}    
\newcommand{\spHeatCap}{c_{\mathrm{p}}}

\newcommand{\TOutside}{\temperature_\infty}
\newcommand{\TOutsideTime}{\temperature_{\infty,\timeVar}}
\newcommand{\dTinf}{\theta}                 
\newcommand{\dTinfTimeIJ}{\theta_{t,ij}}  
\newcommand{\heat}{\dot{Q}}

\newcommand{\length}{L} 
\newcommand{\diameter}{d}
\newcommand{\diameterMin}{\diameter_{\mathrm{min}}}
\newcommand{\diameterDiscrete}{D}
\newcommand{\volume}{V}

\newcommand{\rough}{\epsilon} 
\newcommand{\ratioInsul}{r} 
\newcommand{\condInsul}{\lambda_{\mathrm{i}}} 
\newcommand{\condGround}{\lambda_{\mathrm{g}}} 

\newcommand{\subScriptOuterD}{o}  
\newcommand{\depthPipe}{h} 

\newcommand{\hydrR}{R}  
\newcommand{\frictionfactor}{f} 	
\newcommand{\thermR}{U} 

\newcommand{\inflowID}{a}    
\newcommand{\outflowID}{b}

\newcommand{\valveRhydr}{\zeta}   

\newcommand{\bypValve}{\alpha}     

\newcommand{\lmtd}{LMTD}
\newcommand{\heaterCoef}{\xi}	
\newcommand{\heaterExp}{n} 
\newcommand{\THouse}{\dTinf_{\textrm{house}}}

\newcommand{\Qdemand}{\dot{Q}_{\mathrm{d}}}
\newcommand{\Qdemandi}[1]{\dot{Q}_{\mathrm{d},#1}}
\newcommand{\DemandSatisfaction}{S}

\newcommand{\prodInputi}[1]{\prodInput_{#1}} 	


\newcommand{\existance}{\phi}
\newcommand{\exPipe}{\phi}
\newcommand{\exProducer}{\omega}

\newcommand{\heatInstalled}{H}
\newcommand{\flowVelocity}{v}
\newcommand{\massFlow}{m}
\newcommand{\bigM}{\mathcal{M}}
\newcommand{\conSimpleHX}{c}

\newcommand{\segment}{s}

\newcommand{\gams}{fMINLP}
\newcommand{\pathopt}{pNLP}

\newcommand{\walltime}{w}
\newcommand{\walltimeGAMS}{\walltime_{\mathrm{\gams{}}}}
\newcommand{\walltimePATHOPT}{\walltime_{\mathrm{\pathopt{}}}}

\newcommand{\numberPipes}{n}
\newcommand{\pen}{\xi} 
\newcommand{\penDirection}{a}
\newcommand{\penDiameter}{\bar{\diameter}}
\newcommand{\setDiameters}{\mathcal{S}}
\newcommand{\relu}{f}


\newcommand{\defSetPipes}{\forall \gi\gj \in \Epipe}
\newcommand{\defSetRad}{\forall \gi\gj \in \Erad}
\newcommand{\defSetByp}{\forall \gi\gj \in \Ebyp}
\newcommand{\defSetCon}{\forall \gi\gj \in \Erad \cup \Ebyp}
\newcommand{\defSetProd}{\forall \gi\gj \in \Epro}

\newcommand{\defInflow}{\inflowID=(\gi,n) \in \setEdges}
\newcommand{\defOutflow}{\outflowID=(n,\gj) \in \setEdges}

\newcommand{\topVarBoxConstraints}[1]{\diameterDiscrete_{0}\leq#1\leq\diameterDiscrete_{N}}
\newcommand{\designVarUpper}{\designVar_{\mathrm{up}}}
\newcommand{\designVarLower}{\designVar_{\mathrm{low}}}
\newcommand{\opVarBoxConstraints}[1]{\designVarLower\leq #1 \leq \designVarUpper}
\newcommand{\defFlow}{\ve{\flow} }
\newcommand{\defPressure}{\ve{\pressure} }
\newcommand{\defTemp}{\ve{\dTinf} }

\newcommand{\prodTempSupplyCIJ}{T_{ij,t}}
\newcommand{\prodTempReturnCIJ}{T_{i,t}}

\newcommand{\defDesignVar}{\designVar = \tp{\left[\tp{\ve{\capVar}},\tp{\ve{\radValve}},\tp{\ve{\prodInput}},\tp{\ve{\prodTemp}}\right]}}

\newcommand{\defDesignVarTime}{\designVarTime = \tp{\left[\tp{\ve{\capVar}},\tp{\ve{\radValve}_{\timeVar}},\tp{\ve{\prodInput}_{\timeVar}},\tp{\ve{\prodTemp}_{\timeVar}}\right]}}

\newcommand{\defTopVar}{\topVar\in\{\diameterDiscrete_{0},\dots,\diameterDiscrete_{N}\}^{\card{\Epipe}}}

\newcommand{\defModelConstraints}{\equalCon_\timeVar\left( \designVarTime,\stateVar_\timeVar\right)}

\newcommand{\interestRate}{r_{\mathrm{in}}}
\newcommand{\maxPressure}{\Delta p_\mathrm{max}}
\newcommand{\peakPeriod}{t_{\textrm{peak}}}
\newcommand{\setPeriods}{\Upsilon}

\newcommand{\maxPressureValue}{10}

\clearpage 

\section{Introduction}
\subsection{Background}

District Heating Networks (DHNs) are a network technology that connects heat demand and supply through a network of insulated pipes carrying hot water \cite{WackTopology}. It is considered one of the key technologies for carbon-neutral space heating because of its ability to connect a variety of different renewable and waste heat sources, and provide heat to districts and entire cities \cite{Tilia2022}. However, existing DHNs still rely heavily on fossil fuels. While the pipes of DHNs typically exceed a lifetime of 30+ years \cite{Hay2021}, there is a need to replace or retrofit heat producers in existing networks to enable low-emission and energy efficient heat supply and thus 4th Generation District Heating ($\FGDH$ or $\FG$). 
As low-emission heat sources, renewables and waste heat, are often available at relatively low temperatures, the operating temperatures of DHNs have been steadily decreasing over the last decades. In fact, reducing the supply and return temperature of DHNs is one of the key requirements to enable 4GDH \cite{Volkova2018}. Moreover, lowering the operating temperature of a district heating network (DHN) also reduces its heat losses, which further contributes to the objective of lower network temperatures. However, from a DHN operator's point of view, the operating temperature itself is not necessarily of interest, but rather part of optimizing the overall costs and $\COTWO$ emissions, while ensuring the reliability of the heat supply. 
In this context, determining the optimal network temperature is a complex problem due to multiple reasons: First, when lowering supply temperatures in existing networks, a trade-off between lower heat production costs, increased pumping costs, and exceeding pressure limitations in the network has to be made. Second, today there is also a trade-off to be made between the $\COTWO$ emissions of the heat production units and their CAPEX and OPEX, which is again linked to the network supply temperature and flow rates. This trade-off is particularly challenging given the uncertainty of future $\COTWO$ limits and prices. These complex trade-offs require consideration of heat producer temperatures along with the network temperature, flow, and pressure distributions when retrofitting existing networks for 4G operation. The overall aim is to achieve an optimal balance between costs and $\COTWO$ emissions. This complexity in the design calls for a detailed model-based approach in order to assess the impact of e.g. changing flow temperatures on costs and $\COTWO$ emissions, but also to ensure the feasibility of meeting the heat demand and obeying pressure limits. Variations in heat demand and environmental conditions over the year need to be taken into account to enable the integration of intermittent renewable heat sources and also to achieve more cost-effective network designs by avoiding over-conservatism. For medium to large heating networks, the mentioned characteristics and constraints of the optimal producer design result in a large-scale non-linear optimization problem.

\subsection{Previous work on heat producer retrofit optimization for DHNs}
To solve the problem of selecting and designing new heat producers for an existing DHN (producer retrofit), optimization approaches used in the past have been either parameter optimization, sensitivity studies, or linearized MILP approaches. 
\citet{Popovski2019} compared different, predefined heat producer setups (with fixed capacities), such as solar thermal and heat pumps, but also different renovation and expansion scenarios, for an existing DHN in Germany. The authors used the simulation software energyPRO to evaluate the costs and $\COTWO$ emissions of the different heat supply scenarios. They also conducted a sensitivity study for a scenario that included a heat pump for some of the influencing parameters, such as the price of electricity, to assess how they affect the overall cost. The DHN was not physically modeled but it was instead represented as a single demand point.
\citet{Velasco2023} conducted a simulation-based (TRNSYS) parameter optimization study to compare different, predefined heat producer setups, e.g., a gas boiler + a biomass boiler with and without heat storage. For each setup, the authors assessed different producer/storage capacities with respect to their $\COTWO$ emissions and costs for an existing DHN in Poland. The DHN was not physically modeled but it was instead represented as a single demand point.
\citet{Povzgaj2023} conducted a simulation-based comparison study of three different, predefined producer setups with fixed capacities. The authors goal was to evaluate how the installation of heat pumps for existing DHNs impacts the cost and $\COTWO$ emissions for an existing DHN in Croatia. The authors modeled the DHN based on the simulation software TRNSYS to assess the thermal network losses.
\citet{Lerbinger2023} conducted a MILP optimization to identify optimal decarbonization strategies for two districts in Switzerland, where they investigated the trade-off between $\COTWO$ emissions and costs by creating a Pareto front for each district. An existing DHN was part of the overall decarbonization optimization study for which the sizes of newly build heat producers such as biomass boilers and heat pumps were assessed, next to building renovation strategies and carbon capture technologies. The DHN was not physically modeled but predefined factors were used instead to estimate the total heat loss of the network.
\\
None of the existing work on heat producer retrofit optimization for DHNs demonstrates an automated and physics-based design approach that accurately captures the non-linear physics of heat transfer in the network while applying mathematical optimization to retrofit the heat producers. 
Finding the optimal producer design and operation can only be achieved if the two following requirements are met: First, feasibility with regard to heat demand satisfaction and network pressure limits is ensured by accurately modeling the non-linear heat and pressure losses throughout the DHN while accounting for the existing pipe infrastructure. Second, the capacities and supply temperatures of the potential producers must be optimally determined through mathematical optimization, while accounting for the non-linear temperature dependent producer efficiencies, and not predefined/guessed or only partially explored via a parameter optimization.
\\
As shown in this section, the literature on heat producer retrofit optimization for DHNs is very limited today, that is why we provide an additional literature overview in the next section on the more general heat producer optimization for DHNs, i.e. without the specific requirement of producer retrofit for existing DHNs. This provides an additional reference on used producer models, defined objective functions and optimization approaches. However, note that this is a different problem than producer retrofit for existing DHNs as it does not require the consideration of the existing piping infrastructure.

\subsection{Previous work on heat producer optimization for DHNs (no producer retrofit)}

To solve the problem of selecting and designing heat producers for DHNs, optimization approaches that have been used in the past can be split into two main groups. In the first group, the DHN is either not modeled at all (considered as single demand) or it is modeled but the non-linear nature of the heat transport problem, the temperature dependent efficiencies and the costs have been linearized (simplified DHN modeling approaches). Whereas the second group does model the DHN and account for (some) non-linearities (non-linear DHN modeling approaches).

\subsubsection{Simplified DHN modeling approaches}

\citet{Morvaj2016} solved a multi-objective MILP, considering costs and $\COTWO$ emission limits, for a small DHN to optimally select the heat production technologies, their respective capacities, the network layout and the optimal operation.
\citet{Vesterlund2017} solved a hybrid evolutionary-MILP optimization to determine the optimal supply temperatures and heat inputs of existing heat producers for a large-scale DHN in the city of Kiruna, Sweden. The authors minimized the operating costs of the network.
\citet{Abokersh2020} assessed the possibility of integrating a heat pump into a solar assisted hypothetical DHN in Spain with seasonal thermal energy storage. The authors used a heuristic optimization routine based on a neural network to conduct a multi-objective optimization to balance economical, energetic and environmental goals. The DHN was not modeled in this work.
\citet{Wirtz2021} proposed a MILP for short-term (control) cost optimization of a 5GDHC system, aiming at optimal operation of all generation and storage units of the system as well as optimal network temperature profile.
\citet{Cin2024} presented a MILP for multi-energy system optimization (heat + electricity) of a small district in Italy consisting of eight buildings, aiming at minimizing the life cycle cost while obeying a $\COTWO$ constraint. The authors optimized the additional capacity of already available energy conversion plants and storage units. In one of their scenarios, they considered the development of a DHN.
\citet{Guo2024} proposed a MILP to optimize the DHN size (scope) and its components simultaneously by using spatial clustering techniques. The authors assessed for a district with 255 buildings in Germany how the trade-off between centralized (large) and decentralized (small) DHNs impacts the economic, environmental and technical KPI's and how the chosen DHN vary.
\\ 
The transformation of DHNs towards multi-source, low-temperature networks, breaks the assumptions of most linear heating network models that are used to aid in the producer design. These linearized approaches do not (accurately) account for temperature dependent producer efficiencies or temperature and pressure dependent network distribution losses.
\subsubsection{Non-linear DHN modeling approaches}
\citet{Fang2015} optimized the supply temperatures and heat inputs for a CHP unit and a natural gas boiler to minimize the operating costs. They formulated a non-smooth cost function and used a genetic algorithm to solve a small case consisting of two producers and six consumers.
\citet{Mertz2016} solved a MINLP via GAMS to optimize the network layout, the producer types and their capacities, the heat exchangers and the state variables like temperatures and mass flow rates, aiming at minimizing the overall costs. The authors solved a small case study consisting of four consumers and one producer. Similar to the work of \citet{Mertz2016}, \citet{Marty2018} also conducted a MINLP optimization solved via GAMS but to determine the distribution between electricity and heat production of a geothermal plant. \citet{Maximov2021} conducted a simulation-based multi-objective capacity optimization by using a heuristic genetic algorithm. The assessed systems consisted of solar thermal collectors and thermal storages, with costs and greenhouse gas emissions as objectives for two locations in Chile, while the larger system consisted of 144 buildings.

\subsection{Goal and scope of this paper}
The producer retrofit of existing heating networks with high $\COTWO$ emissions and their transition to $\FGDH$ creates a growing need for scalable, and automated design tools. This paper presents an automated design approach to the optimal producer retrofit of existing networks, based on a detailed physics-based model of future operation, to determine the producer capacities, supply temperatures, and heat loads. The approach accurately captures the spatial and temporal details of the operation, such as temperature and pressure degradation throughout the network, while minimizing costs and $\COTWO$ emissions of the heat producers and pumps. In comparison to previously presented heat producer retrofit approaches, physics-based modeling of the DHN and optimal producer retrofit are ensured through a non-linear model-based mathematical optimization routine. The continuous problem is solved using an adjoint-based Quasi-Newton algorithm. This optimization methodology was already used successfully in previous papers by the authors of this work \cite{Blommaert2020, WackTopology} to optimize the topology of DHNs, and the scalability of this method has been proven \cite{WackBenchmark}. The optimization is formulated as a multi-period problem which allows resolving temporal variations of key properties such as heat demand, outdoor temperature, and solar irradiance \cite{WackMP}. Furthermore, the non-linear physics model allows to retrofit the producers even for DHNs with multiple producers, meshed pipe topologies, and low-temperature operation.
\\
The contributions of this paper are organized as follows. First, the producer retrofit optimization problem for DHNs is defined. Second, producer-specific costs and efficiencies are presented. Third, the methods for efficiently solving this problem are briefly presented, together with the used time aggregation method to manage the computational complexity of the multi-period formulation. Finally, a case study is presented to demonstrate the potential of an automated model-based retrofit design tool for the decarbonization of existing DHNs. Here, we study the impact of a $\COTWO$ price on the producer design, the costs, and the $\COTWO$ emissions in detail.

\section{A multi-objective heat producer retrofit optimization framework}\label{sec:OptProblem}
\newcommand{\gEdge}[3]{#1_{#2#3}}   
\newcommand{\gNode}[2]{#1_{#2}}     

\newcommand{\dirGraph}{G}
\newcommand{\setNodes}{N}
\newcommand{\setEdges}{E}

\newcommand{\pro}{\mathrm{pr}}
\newcommand{\con}{\mathrm{con}}
\newcommand{\rad}{\mathrm{hs}}
\newcommand{\byp}{\mathrm{bp}}
\newcommand{\jun}{\mathrm{jun}}
\newcommand{\pipe}{\mathrm{pipe}}
\newcommand{\operation}{\mathrm{op}}
\newcommand{\Npro}{\setNodes_\pro}
\newcommand{\Ncon}{\setNodes_\con}
\newcommand{\NconF}{\setNodes_{\con,\mathrm{f}}}
\newcommand{\NconR}{\setNodes_{\con,\mathrm{r}}}
\newcommand{\Njun}{\setNodes_\jun}
\newcommand{\NproF}{\setNodes_{\pro,\mathrm{f}}}
\newcommand{\NproR}{\setNodes_{\pro,\mathrm{r}}}

\newcommand{\EF}{\setEdges_{\mathrm{f}}}
\newcommand{\Epro}{\setEdges_\pro}
\newcommand{\EproGB}{\setEdges_{\pro,\GBIndex}}
\newcommand{\EproHP}{\setEdges_{\pro,\HPIndex}}
\newcommand{\EproEB}{\setEdges_{\pro,\EBIndex}}
\newcommand{\EproST}{\setEdges_{\pro,\STIndex}}

\newcommand{\Econ}{\setEdges_\con}
\newcommand{\Erad}{\setEdges_\rad}
\newcommand{\Ebyp}{\setEdges_\byp}
\newcommand{\Epipe}{\setEdges_\pipe}
\newcommand{\EpipeF}{\setEdges_{\pipe,\mathrm{f}}}
\newcommand{\EpipeR}{\setEdges_{\pipe,\mathrm{r}}}
\newcommand{\Eop}{\setEdges_\operation}

\newcommand{\gi}{i}
\newcommand{\gj}{j}

\newcommand{\giNode}[1]{\gNode{#1}{\gi}}
\newcommand{\gjNode}[1]{\gNode{#1}{\gj}}
\newcommand{\gijEdge}[1]{\gEdge{#1}{\gi}{\gj}}

\newcommand{\gjNodetime}[1]{\gNode{#1}{\gj,\timeVar}}
\newcommand{\giNodetime}[1]{\gNode{#1}{\gi,\timeVar}}
\newcommand{\gijEdgetime}[1]{\gEdge{#1}{\gi}{\gj},_\timeVar}

\newcommand{\cost}{J}
\newcommand{\equalCon}{\ve{c}}
\newcommand{\inEqualConSkalar}{\ve{h}}
\newcommand{\inEqualConSkal}{h}
\newcommand{\EqualConSkal}{g}
\newcommand{\inEqualCon}{\inEqualConSkalar}
\newcommand{\designVarskal}{\varphi}
\newcommand{\designVar}{\ve{\designVarskal}}
\newcommand{\designVarTimeInvariantSkalar}{\phi}
\newcommand{\designVarTimeInvariant}{\ve{\designVarTimeInvariantSkalar}}
\newcommand{\designVarTimeDependentSkalar}{\bar{\phi}}
\newcommand{\designVarTimeDependent}{\ve{\designVarTimeDependentSkalar}}
\newcommand{\stateVarSkal}{x}
\newcommand{\stateVar}{\ve{\stateVarSkal}}
\newcommand{\stateVarTime}{\ve{\stateVarSkal}_{\timeVar}}
\newcommand{\topVarSkalar}{d}
\newcommand{\topVar}{\ve{\topVarSkalar}}

\newcommand{\capVar}{\phi} 	
\newcommand{\prodInput}{\gamma} 	
\newcommand{\radValve}{\alpha}    
\newcommand{\prodTemp}{\tau}    
\newcommand{\timeVar}{t}

\newcommand{\n}{n}
\newcommand{\npipes}{\n_{\pipe}}
\newcommand{\nperiods}{\n_{\textrm{period}}}
\newcommand{\timeSlices}{\nperiods}

\newcommand{\defnPeriods}{\nperiods \in \mathbb{N}_1}

\newcommand{\nElements}{\n_{\textrm{spat}}}
\newcommand{\techCon}{\ve{\inEqualConSkal}_{\mathrm{tech},\timeVar}}
\newcommand{\techConEquality}{\ve{\EqualConSkal}_{\mathrm{tech},\timeVar}}
\newcommand{\defStateConstraints}{\techCon(\designVarTime,\stateVarTime)}
\newcommand{\defStateConstraintsEquality}{\techConEquality(\designVarTime,\stateVarTime)}

\newcommand{\nDiscretePipes}{\n_{\mathrm{D}}}
\newcommand{\setDefDiscreteDiameters}{\{\diameterDiscrete_{0},\dots,\diameterDiscrete_{\nDiscretePipes}\}}
\newcommand{\defnpipes}{\npipes = \card{\Epipe}}

\newcommand{\GBName}{\mathrm{gas \ boiler}}
\newcommand{\HPName}{\mathrm{heat \ pump}}
\newcommand{\STName}{\mathrm{solar \ thermal}}
\newcommand{\EBName}{\mathrm{electric \ boiler}}
\newcommand{\GBNameITALIC}{gas \ boiler}
\newcommand{\HPNameITALIC}{heat \ pump}
\newcommand{\STNameITALIC}{solar \ thermal}
\newcommand{\EBNameITALIC}{electric \ boiler}
\newcommand{\GBIndex}{\mathrm{GB}}
\newcommand{\HPIndex}{\mathrm{HP}}
\newcommand{\STIndex}{\mathrm{ST}}
\newcommand{\EBIndex}{\mathrm{EB}}

In this section, the multi-objective heat producer retrofit optimization problem
is formulated and the solution procedure is described. The DHN simulation and optimization is based on previous developments by the authors of this paper \cite{WackTopology, WackMP}. Therefore, the focus of this section is on the new developments of this work, to enable producer retrofit and the multi-objective optimization. For completeness, the main structure of the underlying optimization problem is outlined.

\subsection{Notation}
DHNs are a network technology and can be depicted with a directed graph $\dirGraph(\setNodes,\setEdges)$. Here, $\setNodes$ denotes the set of nodes, while $\setEdges$ represents the set of edges present in the graph. An overview of the graph notation used to formulate optimization problems in this paper is given in table \ref{tab:DHNnotation}.

\begin{table}[h!]
	\centering
	\caption{Graph and index notation for the DHN.}
	\label{tab:DHNnotation}
	\begin{tabularx}{\columnwidth}{XX}
		\toprule
		\textbf{Notation} & \textbf{Description} \\
		\midrule
		$\dirGraph = (\setNodes, \setEdges)$ & DHN as a directed graph.\\
		$\setNodes = \Npro \cup \Ncon \cup \Njun$ & Set of nodes, including heat
		producers, consumers, and
		junctions.\\
		$\setEdges = \Epro \cup \Econ \cup \Epipe$ & Set of edges, including heat
		producers, consumers, and
		pipes.\\
		$\Epro = \EproGB \cup \EproHP \cup \EproST \cup \EproEB$ & Producer subsets, one for each technology.\\
		$\GBIndex$ & $\GBName$\\
		$\HPIndex$ & $\HPName$\\
		$\STIndex$ & $\STName$\\
		$\EBIndex$ & $\EBName$\\
		$(\gi,\gj)$ or $\gi\gj$ & Directed edge going from node $\gi$ to node $\gj$. \\
		\bottomrule
	\end{tabularx}
\end{table}

\subsection{Heat producer retrofit as an optimization problem}
First, the multi-objective optimization problem for heat producer retrofit is formulated. Designing new producers for a DHN requires to consider the performance throughout future operation. Therefore, optimizing the producer retrofit means solving a multi-period optimization problem to represent multiple operating points throughout the year. To account for the time-dependent nature of this optimization problem, a set of time periods $\setPeriods=\{1,\dots,\timeSlices\}$ is defined, where $\defnPeriods$ is the number of time periods. An individual period is defined by the index $\timeVar \in \setPeriods$. The worst-case period during which the network must ensure feasible operation is denoted by $\peakPeriod\in\mathbb{N}_1$. 
\\
In the corresponding optimization problem, the project cost $\costFull$ is minimized by choosing for each considered producer type the capacity $\ve{\capVar}$, and for each period the producer supply temperature $\ve{\prodTemp}_{\timeVar}$ and the inflow $\ve{\prodInput}_{\timeVar}$. To ensure a feasible operation, the valve settings of the heat consumer substations $\ve{\radValve}_{\timeVar}$ in each time period $\timeVar$ have to be chosen. The multi-period producer retrofit optimization problem for an existing DHN is thus

\begin{equation}
	\begin{aligned} \label{eq:optMP}
		\min_{\designVar,\stateVar} &\qquad
		\costFull \left(\designVar,\stateVar\right) &\\
		s.t.& \qquad \defModelConstraints = 0,\quad\forall \timeVar \in \setPeriods\,, &\\
		& \qquad \defStateConstraints \leq 0,\quad\forall \timeVar \in \setPeriods\,, & \\
		& \qquad \defStateConstraintsEquality = 0,\quad\forall \timeVar \in \setPeriods\,, & \\
		& \qquad \opVarBoxConstraints{\designVar}\,, &
	\end{aligned}
\end{equation}

\newcommand{\defStateVar}{\ve{\stateVar}_{\timeVar} = \tp{\left[\tp{\ve{\flow}}_{\timeVar},\tp{\ve{\pressure}}_{\timeVar},\tp{\ve{\dTinf}}_{\timeVar}\right]}} 

where the design variables $\defDesignVarTime$ and the physical variables $\defStateVar$ are optimized to satisfy the model equations $\defModelConstraints$. For ease of notation, the time-dependent variable vectors $\designVarTime$ and $\stateVar_\timeVar$ of all periods are combined in the vectors $\designVar$ and $\stateVar$, respectively. The set of non-linear model equations represents the hydraulic and thermal transport problem in the network. It includes the flow rates $\defFlow_{\timeVar}$, nodal pressures $\defPressure_{\timeVar}$, and nodal and pipe exit temperatures $\defTemp_{\timeVar}$.
The temperature $\defTemp_{\timeVar} = \ve{\temperature}_{\timeVar} - {\temperature}_{\infty,\timeVar}$ is defined as the difference between the absolute water temperature $\ve{\temperature}_{\timeVar}$ and the outside air temperature ${\temperature}_{\infty,\timeVar}$. Additional technological constraints (e.g. meeting consumer heat demands), are represented by inequality $\defStateConstraints$ and equality $\defStateConstraintsEquality$ state constraints.

\newcommand{\npv}{NPV}
\newcommand{\npvConst}{f}
\newcommand{\npvN}{A}
\newcommand{\npvt}{k}
\newcommand{\npvDiscount}{e}

\newcommand{\npvConstCAPEX}{\npvConst_{\mathrm{\subCAPEX}}}
\newcommand{\npvConstOPEX}{\npvConst_{\mathrm{\subOPEX}}}

\newcommand{\npvCost}{C}

\newcommand{\actuarialRate}{\npvDiscount_{\mathrm{a}}}
\newcommand{\energyInflation}{\npvDiscount_{\mathrm{i}}}

\newcommand{\periodWeight}{\omega}

The discounted lifetime cost of the producer retrofit $\costFull$ in this multi-objective approach includes both the CAPEX and OPEX, and the $\COTWO$ costs of the heat producers ($\costi{\subHeat}$) and pumps ($\costi{\subPump}$), and is defined as

\begin{align}
	\costFull\left(\designVar,\stateVar\right) &=  \sum_{\gi\gj \in \Epro}\costi{\subHeat,\subCAPEX, \textit{ij}} \nonumber\\
	&+ \npvConstOPEX \sum_{\timeVar = 1}^{\timeSlices} \periodWeight_{\timeVar} \sum_{\gi\gj \in \Epro} \left[\costi{\subHeat,\subOPEX, \textit{ij},\textit{\timeVar}} + \costi{\subHeat, \subCOTWO, \textit{ij},\textit{\timeVar}}\right]
	\nonumber\\
	&+ \npvConstOPEX \sum_{\timeVar = 1}^{\timeSlices} \periodWeight_{\timeVar} \sum_{\gi\gj \in \Epro} \left[\costi{\subPump,\subOPEX, \textit{ij},\textit{\timeVar}} + \costi{\subPump, \subCOTWO, \textit{ij},\textit{\timeVar}}\right]
	\label{eq:totalCost}\,,
\end{align}

with $\npvConstOPEX = \sum_{\npvt=1}^{\npvN} \frac{1}{\left(1+\npvDiscount\right)^\npvt} $ assuming an investment horizon of $\npvN = 30~\mathrm{years}$ and a discount rate of $\npvDiscount=5\%$. Here, it is assumed that the chosen time slices $\timeVar$ are sufficiently representative for a full year operation and $\sum_{\timeVar=1}^{\timeSlices}\periodWeight_{\timeVar} = 1$. The individual cost components $\costi{}$ are described in the next sections for the four heat producer technologies considered in this work: $\GBName$, $\HPName$, $\STName$ and $\EBName$.

\newcommand{\cHeatCAPEX}{\npvCost_{\mathrm{hC}}}
\newcommand{\cHeatCAPEXi}[1]{\npvCost_{\mathrm{hC},#1}}
\newcommand{\capacityFactor}{F}
\newcommand{\efficiency}{\eta}
\newcommand{\producerEfficiency}{\efficiency_{\pro}}

\newcommand{\maxCapacity}{P_{\textrm{max}}}
\newcommand{\maxCapacityIJ}{P_{\textrm{max},ij}}
\newcommand{\maxSurfaceST}{A_{\textrm{max}}}
\newcommand{\prodCostSpecific}{\npvCost_{\textrm{C,prod}}}
\newcommand{\prodCostSpecificIJ}{\npvCost_{\textrm{C,prod,\textit{ij}}}}
\newcommand{\prodCostSpecificIJbold}{\bm{\npvCost}_{\textbf{\textrm{C,prod,\textit{ij}}}}}

\newcommand{\etaProd}{\eta_{\textrm{prod}}}
\newcommand{\etaProdTime}{\eta_{\textrm{prod},\timeVar}}
\newcommand{\etaProdIJ}{\eta_{\textrm{prod},ij}}
\newcommand{\etaProdIJTime}{\eta_{\textrm{prod},ij,t}}
\newcommand{\etaProdIJTimebold}{\bm{\eta}_{\textbf{\textrm{prod}}\bm{,ij,t}}}

\subsection{Cost function}

\subsubsection{Heat producer CAPEX}

\begin{figure}[h]
	\includegraphics[width=\linewidth]{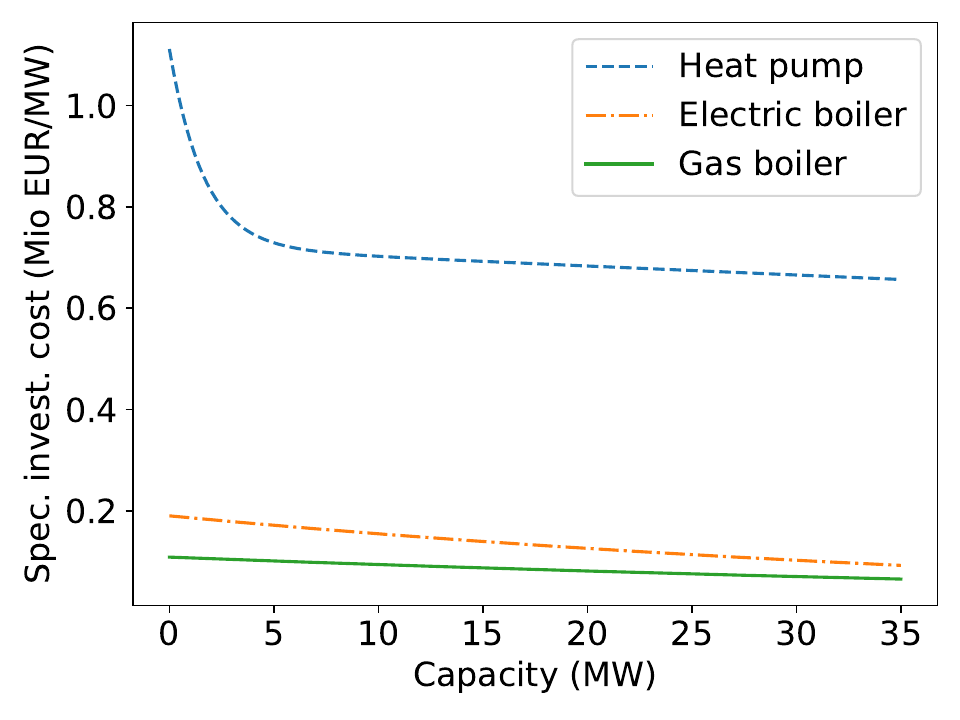}
	\caption{Fitted specific (spec.) investment cost $\prodCostSpecificIJ$ of the $\GBName$, $\HPName$ and $\EBName$ heat producers.}
	\label{fig:specificInvestmentCosts}
\end{figure}

The investment costs (CAPEX) of the heat production components $\GBName$, $\HPName$, and $\EBName$ can be formulated as the product of the specific investment costs $\prodCostSpecific(\capVar\,\maxCapacity)$ (unit $\si{\sieuro}$/W) times the installed capacity $\capVar\,\maxCapacity$ (unit W) as

\begin{multline}\label{eq:capacityCost}
	\costi{\subHeat,\subCAPEX, \textit{ij}} = \gijEdge{\big(\capVar \, \maxCapacity  \, \prodCostSpecific(\capVar \, \maxCapacity)\big)} \\ \quad \forall ij \in \Epro \setminus \EproST, 
\end{multline}

with $\maxCapacity$ being the maximal capacity. Note that the fitted specific investment costs already incorporate the fixed investment costs. The  specific investment costs are usually decreasing with increasing capacity which is also known as the economies of scale. Moreover, this decrease of the specific investment costs is often non-linear which further highlights the importance of a non-linear and thereby accurate optimization framework. For this work, academic literature and manufacturer data was used to fit the specific investment costs of all heat producer technologies. Figure \ref{fig:specificInvestmentCosts} shows the specific investment costs for the natural gas boiler, air-source heat pump, and electric boiler. The fitted functions and cost references are given in \ref{app:costs}.

\begin{figure}[h]
	\includegraphics[width=\linewidth]{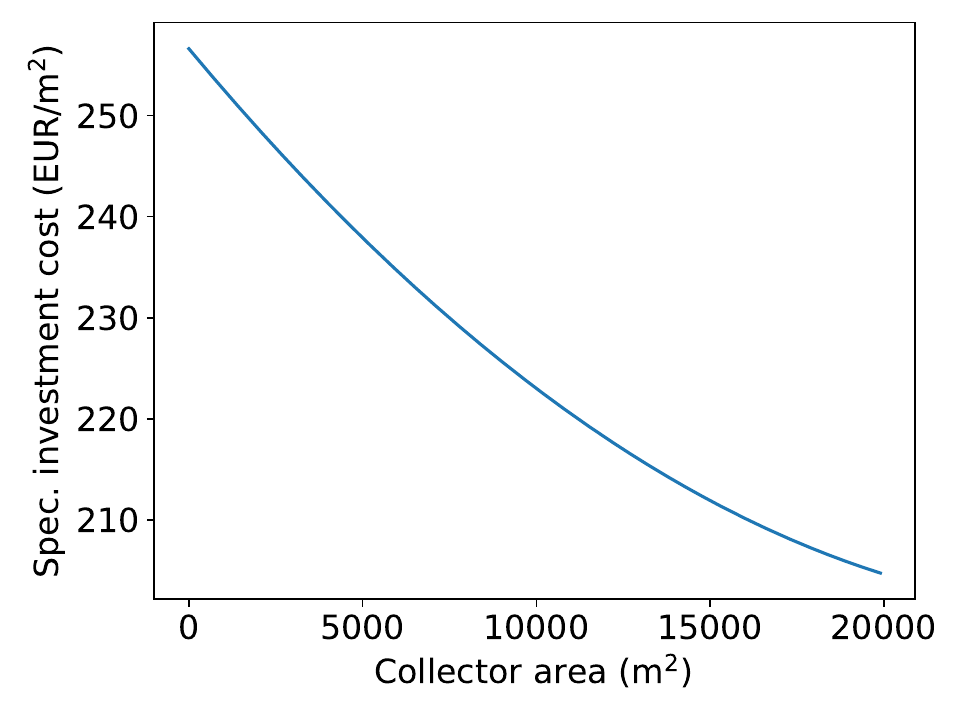}
	\caption{Fitted specific (spec.) investment cost $\prodCostSpecificIJ$ of the heat producer $\STName$.}
	\label{fig:specificInvestmentCostSolar}
\end{figure}

For $\STName$ collectors, the specific investment costs are usually not given per installed peak capacity since the heat output strongly depends on the location and the operation. Instead, they are given per square meter of collector area. Moreover, since large areas of land are needed for the installation of $\STName$ collectors meant to supply DHNs, the cost of land must be considered as well. The fitted specific investment costs can be seen in figure \ref{fig:specificInvestmentCostSolar} and the fitted function and cost references are given in \ref{app:costs}. The total investment cost of the $\STName$ unit are given by

\begin{multline}\label{eq:capacityCostST}
	\costi{\subHeat,\subCAPEX, \textit{ij}} = \gijEdge{\big(\capVar \, \maxSurfaceST  \, \prodCostSpecific(\capVar \, \maxSurfaceST)\big)} \\ \quad \forall ij \in \EproST,
\end{multline}

with $\maxSurfaceST$ being the maximal collector area.

\subsubsection{Heat producer OPEX} 
\label{sec:OPEX}
\subsubsection*{Gas boiler, $\HPNameITALIC$ and $\EBNameITALIC$} 
The considered operational costs (OPEX) of the heat producers are the costs of used electricity and natural gas (energy input). For the $\GBName$, $\HPName$, and $\EBName$ the efficiency $\etaProdIJ$ is used to determine the energy input based on the supplied heat to the DHN and to price it with the specific energy prices which reads 

\begin{multline}\label{eq:Produceropex}
	\costi{\subHeat,\subOPEX, \textit{ij}, \textit{\timeVar}} = \density \, \spHeatCap \, \KOPEX  \left( \cHeatOPEXi \, \flow_{\timeVar} \, \Delta\theta_{\timeVar} \, \frac{1}{\etaProdTime} \right)_{ij} \\ \quad \forall ij \in \Epro \setminus \EproST .
\end{multline}

The temperature increase over the producer is given by $\Delta\dTinfTimeIJ$ while the prices $\cHeatOPEXiIJ$ are given in table \ref{tab:EnergyPrices}. The water density is denoted with $\density$ and has a value of 983$\,\mathrm{kg/m^3}$. The specific heat capacity of water is denoted with $\spHeatCap$ and has a value of 4185$\,\mathrm{J/kgK}$. The conversion factor $\KOPEX$ in $\unit{\hour\per\year}$ defines the network's number of active operating hours per year. 

\begin{table}[h!]
	\centering
	\caption{The prices $\cHeatOPEXiIJ$ in $\si{\sieuro}$/kWh for the operational costs of the different heat producer types. The electricity \cite{ElectricityPrice} and natural gas \cite{GasPrice} price are the average prices in the EU from 2019 until mid 2022 for non-household consumers.}
	\label{tab:EnergyPrices}
	\begin{tabularx}{\columnwidth}{>{\hsize=.27\hsize}X>{\hsize=.21\hsize}X>{\hsize=.51\hsize}X}
		\toprule
		\textbf{Input} & \textbf{$\cHeatOPEXiIJbold$} & \textbf{Producer} \\
		\midrule
		Natural gas & 0.0319 & $\forall ij \in \EproGB$\\ 
		Electricity & 0.1 & $\forall ij \in \EproHP \cup \EproEB \cup \EproST$\\
		\bottomrule
	\end{tabularx}
\end{table}

\begin{figure}[h]
	\centering
	\includegraphics[width=1.0\columnwidth]{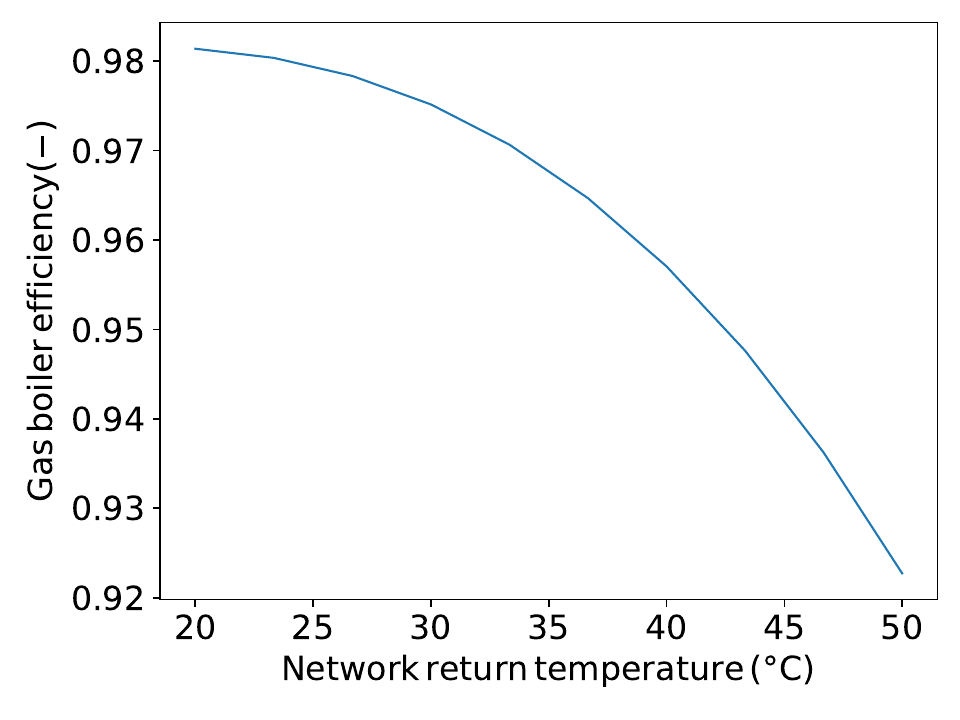}
	\caption{Efficiency curve of the heat producer $\GBName$, showing its dependence on the network return temperature.}
	\label{fig:GB_eta}
\end{figure}

\begin{figure}[h]
	\centering
	\includegraphics[width=1.0\columnwidth]{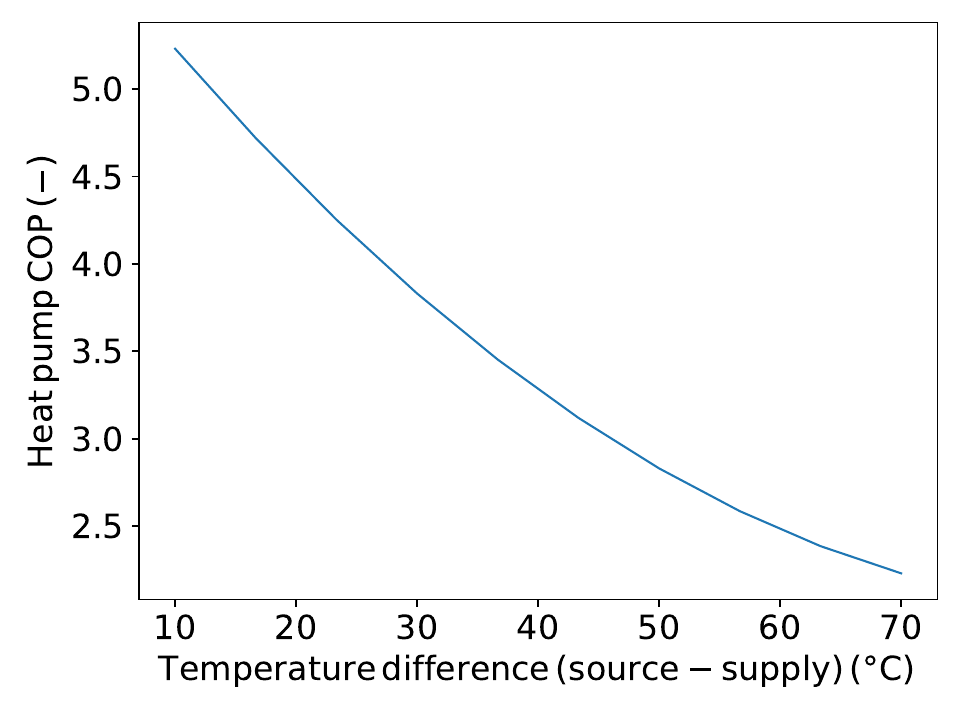}
	\caption{Efficiency (COP) curve of the heat producer $\HPName$, showing its dependence on the temperature difference between the producer supply and the source (air) temperature.}
	\label{fig:HP_eta}
\end{figure}

Efficiencies of heat production units ($\etaProdIJ$) often depend on the operational temperatures, which must be considered to accurately assess the impact of the network temperatures on the costs and $\COTWO$ emissions. For a condensing $\GBName$ the efficiency depends primarily on the inlet temperature (the return temperature from the DHN). Here, values from \citet{Qin2023} are used to fit a function for the efficiency, visualized in figure \ref{fig:GB_eta}. For the air-source $\HPName$ the coefficient of performance (COP) is modeled based on \citet{Ruhnau2019} and depends on the difference between the heat supply and the heat source temperature (see figure \ref{fig:HP_eta}). Concerning the $\EBName$, a fixed efficiency based on \citet{Johansson2020} is assumed. The producer supply temperature in $^\circ\mathrm{C}$ is defined as $\prodTempSupplyCIJ$ which links to the temperature state variable $\theta_{ij,\timeVar}$ as

\begin{equation}\label{eq:producerSupplyTemp}
	\prodTempSupplyCIJ =\theta_{ij,\timeVar} + \TOutsideTime \quad \forall ij \in \Epro,
\end{equation}

whereas the network return temperature (producer inlet temperature) in $^\circ \mathrm{C}$ is given by $\prodTempReturnCIJ$ and its relation to the state variable $\theta_{i,\timeVar}$ is defined as 

\begin{equation}\label{eq:producerReturnTemp}
	\prodTempReturnCIJ =\theta_{i,\timeVar} + \TOutsideTime \quad \forall ij \in \Epro.
\end{equation}

The underlying functions of the efficiency curves are given in table \ref{tab:ProducerEfficiencies}. 

\newcommand{\aGBeta}{a_{\eta,\GBIndex}}
\newcommand{\bGBeta}{b_{\eta,\GBIndex}}
\newcommand{\cGBeta}{c_{\eta,\GBIndex}}
\newcommand{\dGBeta}{d_{\eta,\GBIndex}}
\newcommand{\aHPeta}{a_{\eta,\HPIndex}}
\newcommand{\bHPeta}{b_{\eta,\HPIndex}}
\newcommand{\cHPeta}{c_{\eta,\HPIndex}}
\newcommand{\aEBeta}{a_{\eta,\EBIndex}}

\newcommand{\aGBetavalue}{-7.2225e-07}
\newcommand{\bGBetavalue}{5.5968e-06}
\newcommand{\cGBetavalue}{0.0005}
\newcommand{\dGBetavalue}{0.9755}
\newcommand{\aHPetavalue}{0.0005}
\newcommand{\bHPetavalue}{- 0.09}
\newcommand{\cHPetavalue}{6.08}
\newcommand{\aEBetavalue}{0.98}

\newcommand{\secondaryIndex}{\mathrm{sec}}
\newcommand{\deltaPSEC}{\Delta\pressure_{\secondaryIndex}}
\newcommand{\deltaPSECIJ}{\Delta\pressure_{\secondaryIndex,ij}}

\FloatBarrier
\subsubsection*{Solar thermal} 
For the $\STName$ unit, the OPEX cost differs from those of the $\GBName$, $\HPName$, and $\EBName$ since its main energy input, the solar irradiance, has no cost. However, pumping the heat-transfer fluid through the collector requires electricity which has operational cost. Since the pressure within the collector is not modeled in this work, as will be elaborated in section \ref{section:STModel}, we assume a fixed pressure drop of $\deltaPSECIJ$ = 1 bar. The subscript "$\secondaryIndex$" stands for the secondary side within the heat producer unit (and not within the DHN). The OPEX of the $\STName$ unit is then given by

\begin{multline}\label{eq:STProduceropex}
	\costi{\subHeat,\subOPEX, \textit{ij}, \textit{\timeVar}} = \frac{\KOPEX}{\pumpEff} \, \left(\deltaPSEC \, \flow_{\secondaryIndex,\timeVar} \, \cHeatOPEXi \right)_{ij} \\ \quad \forall ij \in \EproST 
\end{multline}

with $\pumpEff$ being the efficiency of the distribution pump with a value of 0.81.

\subsubsection{Pumping costs}
\label{sec:Pumping}
The operational cost of the distribution pumps at the heat producers is computed with

\begin{multline}\label{eq:Pumpopex}
	\costi{\subPump,\subOPEX, \textit{ij}, \textit{\timeVar}} = \frac{\KOPEX}{\pumpEff} \cPumpOPEXi \, \flow_{ij,\timeVar} \, (\pressure_{j,\timeVar}-\pressure_{i,\timeVar}) \\ \quad \forall ij \in \Epro ,
\end{multline}
while $\cPumpOPEXi$ is the electricity price.

\newcommand{\COTWOcon}{\textit{EF}_{\mathrm{CO_2}}}
\newcommand{\COTWOconIJ}{\textit{EF}_{\mathrm{CO_2},ij}}
\newcommand{\COTWOconIJbold}{\textbf{\textit{EF}}_{\bm{\mathrm{CO_2},ij}}}
\newcommand{\COTWOPrice}{\npvCost_{\mathrm{CO_2}}}

\subsubsection{Pricing $\COTWOITALIC$ emission costs}
The $\COTWO$ costs are directly linked to the previously mentioned OPEX (\ref{sec:OPEX}) and pumping costs (\ref{sec:Pumping}) but instead of pricing the purchase of the energy inputs, natural gas, and electricity, their related $\COTWO$ emissions are priced with a $\COTWO$ price. Hence, the $\COTWO$ costs of the heat producers $\costi{\subHeat}$ and the pumps $\costi{\subPump}$ can be stated as

\begin{multline}\label{eq:ProducerCO2}
	\costi{\subHeat, \subCOTWO, \textit{ij}, \textit{\timeVar}} = \density \, \spHeatCap \, \KOPEX \, \COTWOPrice \left(\COTWOcon \, \flow_{\timeVar} \, \Delta\theta_{\timeVar} \, \frac{1}{\etaProdTime} \right)_{ij} \\ \quad \forall ij \in \Epro \setminus \EproST ,
\end{multline}

\begin{multline}\label{eq:ProducerSTCO2}
	\costi{\subHeat, \subCOTWO, \textit{ij}, \textit{\timeVar}} = \frac{\KOPEX}{\pumpEff} \, \COTWOPrice \left(\COTWOcon \, \deltaPSEC \, \flow_{\secondaryIndex,\timeVar} \, \right)_{ij} \\ \quad \forall ij \in \EproST 
\end{multline}

and,

\begin{multline}\label{eq:PumpCO2}
	\costi{\subPump, \subCOTWO, \textit{ij}, \textit{\timeVar}} = \frac{\KOPEX}{\pumpEff} \, \COTWOPrice \,\COTWOcon \, \flow_{ij,\timeVar} \, (\pressure_{j,\timeVar}-\pressure_{i,\timeVar}) \\ \quad \forall ij \in \Epro .
\end{multline}

The $\COTWO$ price in $\si{\,\sieuro}$/kg$\COTWO$ is given by $\COTWOPrice$ and the specific $\COTWO$ emission factor in kg$\COTWO$/kWh is given by $\COTWOcon$. This factor is dependent on the energy source and the respective values for natural gas and electricity can be found in table \ref{tab:CO2emissions}. For the circulation pumps in equation \ref{eq:PumpCO2} it is always that of electricity.

\begin{table}[h!]
	\centering
	\caption{The specific $\COTWO$ emissions $\COTWOconIJ$ in kg$\COTWO$/kWh of all heat producers. The $\COTWO$ emission value of electricity represents that of Belgian electricity production in the year 2020 \cite{ElectricityCO2}. The $\COTWO$ emission value of natural gas combustion is taken from \cite{GasCO2}.}
	\label{tab:CO2emissions}
	\begin{tabularx}{\columnwidth}{>{\hsize=.25\hsize}X>{\hsize=.18\hsize}X>{\hsize=.42\hsize}X}
		\toprule
		\textbf{Input} & $\COTWOconIJbold$ & \textbf{Producer} \\
		\midrule
		Natural gas & 0.181 & $\forall ij \in \EproGB$\\ 
		Electricity & 0.181 & $\forall ij \in \EproHP \cup \EproEB \cup \EproST$\\
		\bottomrule
	\end{tabularx}
\end{table}

\newcommand{\etaZeroST}{\eta_{0,\STIndex}}
\newcommand{\Girr}{G_{\mathrm{irr}}}
\newcommand{\TmeanST}{T_{m}}
\newcommand{\TreducedST}{T_{m}^*}
\newcommand{\TreducedSTsquare}{T_{m}^{*2}}
\newcommand{\GirrTime}{G_{\mathrm{irr},\timeVar}}
\newcommand{\TmeanSTTime}{T_{\mathrm{m},\timeVar}}
\newcommand{\TreducedSTTime}{T_{\mathrm{m},\timeVar}^*}
\newcommand{\TreducedSTsquareTime}{T_{\mathrm{m},\timeVar}^{*2}}

\subsection{A $\STNameITALIC$ model}
\label{section:STModel}
In this section, we describe the model of the $\STName$ unit. 
\begin{figure}[h!]
	\centering
	\includegraphics[width=\columnwidth]{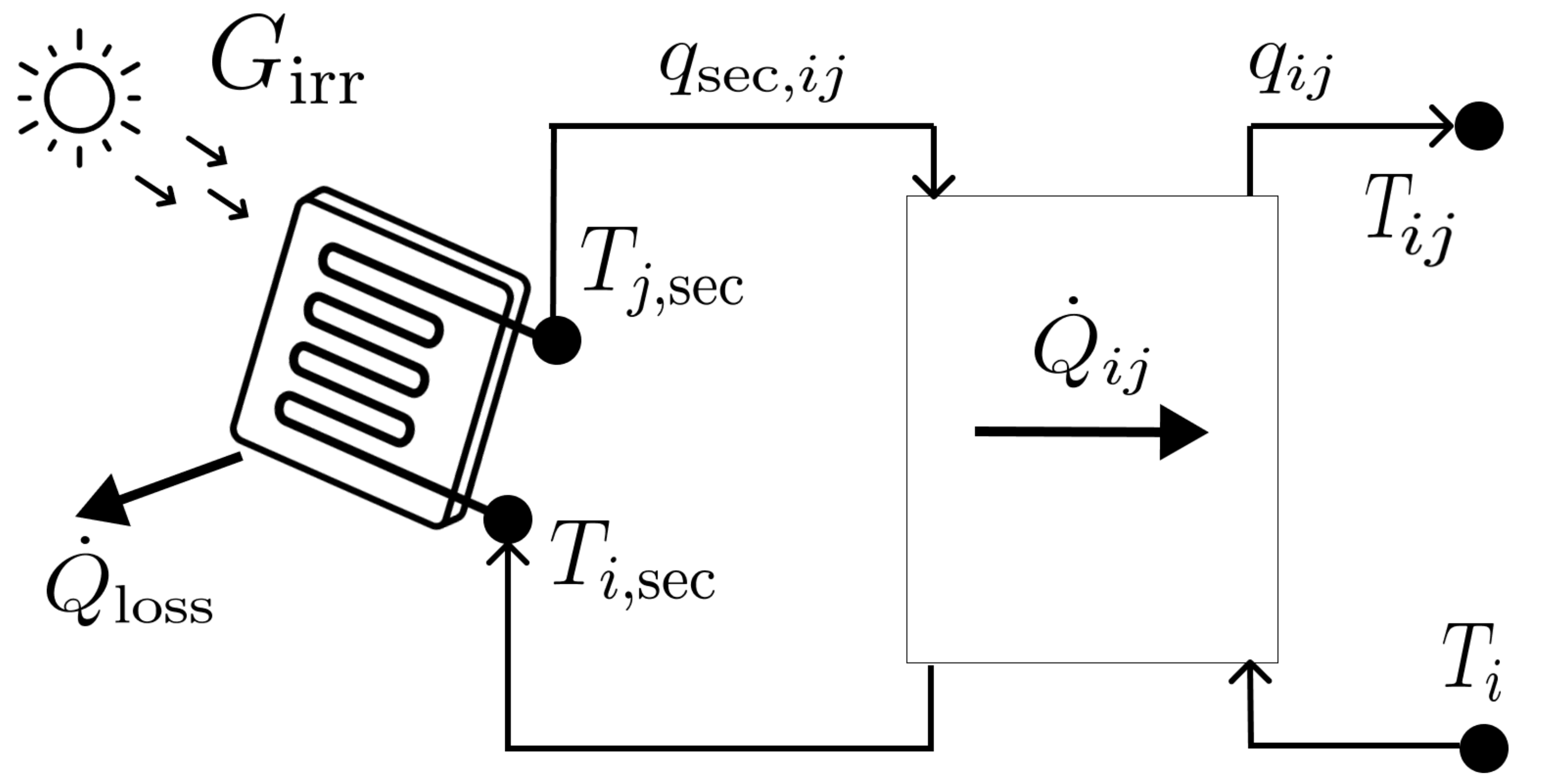} 
	\caption{Visualization of the $\STName$ model and its thermal connection to the DHN via a heat exchanger. Hydraulically, the $\STName$ unit and the DHN are decoupled.}
	\label{fig:solarunit}
\end{figure}
\subsubsection{Background}

The model aims to represent a non-concentrating and insulated flat-plate $\STName$ collector. This collector type is well suited for $\FGDH$ networks because of their supply temperatures which are usually between 40°C - 80°C \cite{Dharuman2006, Hakkarainen2015}. The proposed model is based on real-world measurement data of existing $\STName$ collectors that supply heat to DHNs. The data comes from solar DHN sites in Denmark and is openly accessible \cite{DanishSolarData}. Our model accounts for thermal losses via a global efficiency curve that includes the collector efficiency and the heat exchange between the $\STName$ collector loop and the DHN, since they are hydraulically separated, to reflect existing installations. Moreover, as we are not considering heat storage in this work, we assume that the produced heat of the collector field of each period must be fully integrated into the network to avoid overheating or the need for cooling installations. A visualization of the $\STName$ unit and its connection to the DHN is given in figure \ref{fig:solarunit}.
\subsubsection{Model equations}
\label{subsec:ST_model}
To account for thermal losses of a solar collector its steady-state efficiency curve is commonly used which is also defined by the norm ISO 9806:2014 \cite{SolarNorm}. This curve links the steady-state efficiency
of the collector $\etaZeroST$ to the total solar irradiation $\GirrTime$, ambient
temperature $\TOutsideTime$ and mean solar collector temperature $\TmeanSTTime$. For the efficiency curve, the reduced temperature $\TreducedSTTime = \frac{\TmeanSTTime-\TOutsideTime}{\GirrTime}$ is defined, the efficiency then reads

\newcommand{\aoneSTeta}{a_{1,\STIndex}}
\newcommand{\atwoSTeta}{a_{2,\STIndex}}

\newcommand{\etaZeroSTvalue}{0.75}
\newcommand{\aoneSTetavalue}{3.5}
\newcommand{\atwoSTetavalue}{0.012}

\newcommand{\aTjST}{a_{T_j,\STIndex}}
\newcommand{\bTjST}{b_{T_j,\STIndex}}
\newcommand{\cTjST}{c_{T_j,\STIndex}}
\newcommand{\dTjST}{d_{T_j,\STIndex}}
\newcommand{\eTjST}{e_{T_j,\STIndex}}
\newcommand{\fTjST}{f_{T_j,\STIndex}}
\newcommand{\aTiST}{a_{T_i,\STIndex}}
\newcommand{\bTiST}{b_{T_i,\STIndex}}
\newcommand{\cTiST}{c_{T_i,\STIndex}}
\newcommand{\dTiST}{d_{T_i,\STIndex}}
\newcommand{\eTiST}{e_{T_i,\STIndex}}
\newcommand{\fTiST}{f_{T_i,\STIndex}}

\newcommand{\aTjSTvalue}{23.4137}
\newcommand{\bTjSTvalue}{0.0501}
\newcommand{\cTjSTvalue}{0.7843}
\newcommand{\dTjSTvalue}{1.5235e-05}
\newcommand{\eTjSTvalue}{-2.3705e-04}
\newcommand{\fTjSTvalue}{0.0109}
\newcommand{\aTiSTvalue}{20.9417}
\newcommand{\bTiSTvalue}{0.0216}
\newcommand{\cTiSTvalue}{0.5447}
\newcommand{\dTiSTvalue}{-2.374e-06}
\newcommand{\eTiSTvalue}{-6.6886e-04}
\newcommand{\fTiSTvalue}{0.0173}

\begin{multline}\label{eq:etaST}
	\etaProdIJTime = \etaZeroST - \aoneSTeta\TreducedSTTime - \atwoSTeta\GirrTime\TreducedSTsquareTime \\ 
	\quad \forall ij \in \EproST .
\end{multline}

\begin{figure}[h!]
	\centering
	\includegraphics[width=\columnwidth]{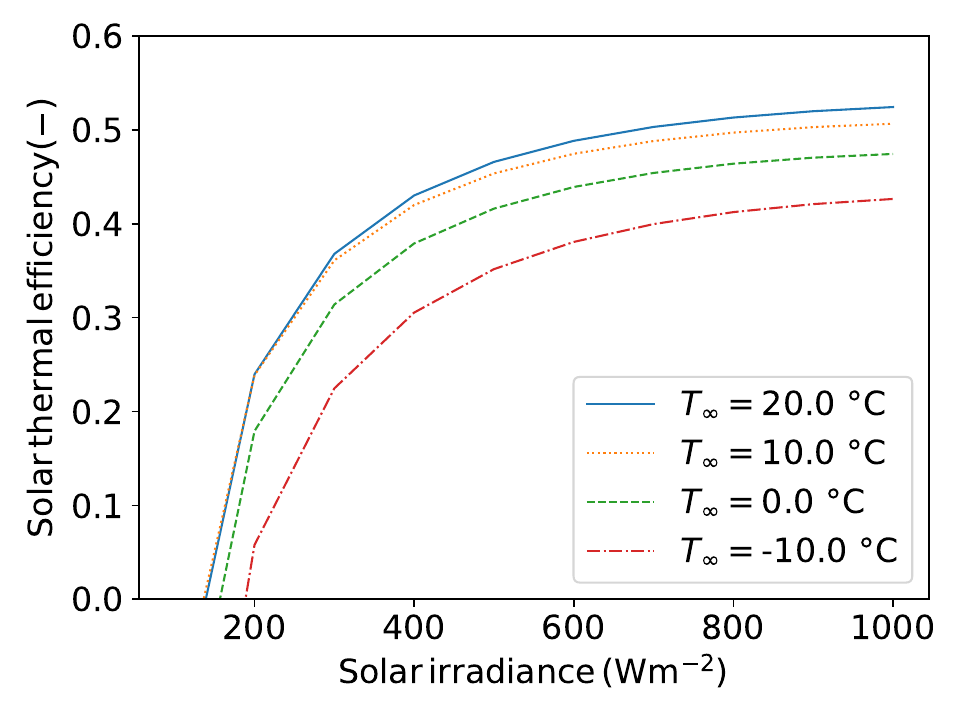}
	\caption{The dependence of the $\STName$ efficiency on the solar irradiance for different outside temperatures.}
	\label{fig:ST_eta}
\end{figure}

\begin{figure}[h!]
	\centering
	\includegraphics[width=\columnwidth]{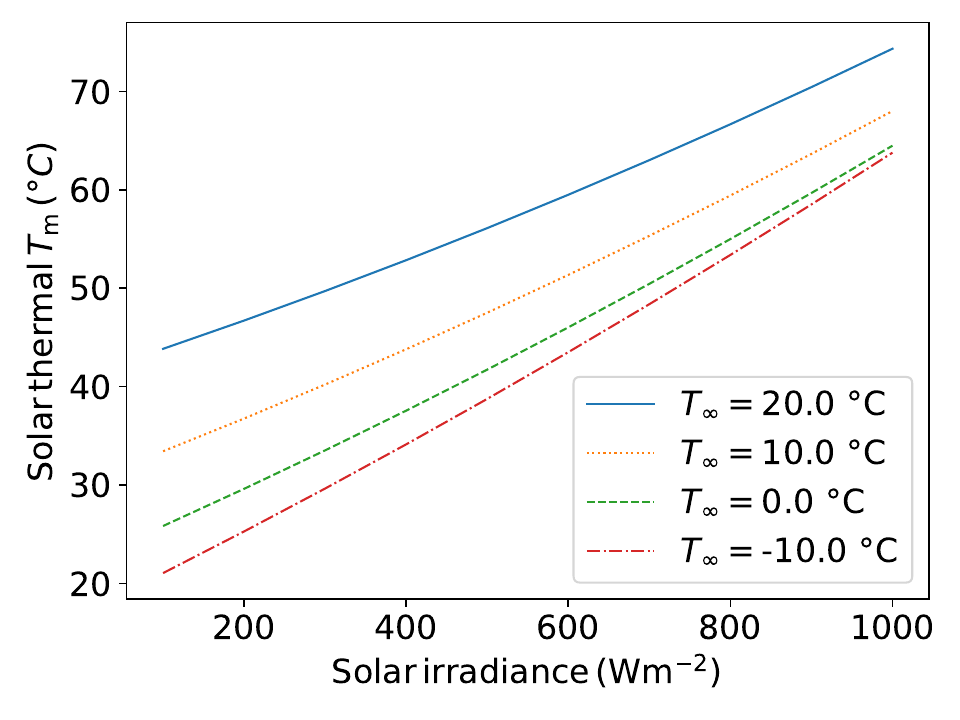}
	\caption{The dependence of the $\STName$ average collector temperature on the solar irradiance for different outside temperatures.}
	\label{fig:ST_Tm}
\end{figure}

The parameters $\etaZeroST$, $\aoneSTeta$ and $\atwoSTeta$ were taken from product data sheets of a solar collector manufacturing company \cite{Gasokol} and are given in table \ref{tab:etaSTParameters}. To obtain the mean solar collector temperature $\TmeanSTTime = \frac{\temperature_{\secondaryIndex,j,\timeVar}+\temperature_{\secondaryIndex,i,\timeVar}}{2}$ we defined fits for the hot $\temperature_{\secondaryIndex,j,\timeVar}$ and the cold $\temperature_{\secondaryIndex,i,\timeVar}$ collector temperature based on the measurement data of existing solar DHNs. These fits are functions of the solar irradiance $\GirrTime$ and the outside temperature $\TOutsideTime$ and allow to realistically estimate the operation of the solar collector field for the different periods of the year. The fits governing equations are given in table \ref{tab:STtemperatureFits} and the related parameter values are given in table \ref{tab:STtemperatureFitsParameters}. Figures \ref{fig:ST_eta} and \ref{fig:ST_Tm} visualize the dependency of the $\STName$ efficiency and mean collector temperature on the solar irradiance and outside temperature. The temperatures and efficiency of the solar collector are pre-processed for each period based on its solar irradiance and outside temperature.

For modeling the heat exchange between the hydraulically separated $\STName$ collector loop and the DHN we assume a counter-flow heat exchanger and use the epsilon-NTU method, similar to the way we model the consumer substations \cite{WackMP}. We assume a 60\% water 40\% propylene-glycol mixture as the heat transfer fluid in the solar collector with a density of 1009 $\mathrm{kg/m^3}$ and a specific heat capacity of 3511 $\mathrm{J/kg\,K}$ based on \cite{Gasokol, Pekasol}. The flow rate is a priori set for each period to allow the full heat integration based on the solar irradiance, the $\STName$ unit's temperatures, and its efficiency. For periods where the solar irradiance and/or the outside temperature are too low to effectively produce heat we set the heat output a priori to zero. The heat exchange is then eventually described by the following three equations

\newcommand{\Qproin}{\dot{Q}_{ij,\timeVar}}
\newcommand{\Qprohx}{\dot{Q}_{\mathrm{hx},ij,\timeVar}}
\newcommand{\Qprosol}{\dot{Q}_{\mathrm{solar},ij,\timeVar}}
\newcommand{\Cmin}{C_{\text{min}}}										
\newcommand{\Cmax}{C_{\text{max}}}										
\newcommand{\Cstar}{C^*}										
\newcommand{\NTU}{NTU}											
\newcommand{\U}{U}												
\newcommand{\epsilonNTU}{\epsilon}								
\newcommand{\CminTime}{C_{\text{min},\timeVar}}										
\newcommand{\epsilonNTUTime}{\epsilon_{\timeVar}}								

\begin{multline}\label{eq:HXST1}
	\Qproin = \density \, \spHeatCap \left(\flow_{\timeVar} \, \Delta\theta_{\timeVar} \, \right)_{ij} \quad \forall ij \in \Epro ,
\end{multline}

\begin{multline}\label{eq:HXST2}
	\Qprohx = \left(\epsilonNTUTime \, \CminTime \, (\temperature_{\secondaryIndex,j,t} - \prodTempReturnCIJ) \right)_{ij} \\ \quad \forall ij \in \EproST ,
\end{multline}

and

\begin{multline}\label{eq:HXST3}
	\Qprosol = \GirrTime \left(\etaProdTime \, \capVar \, \maxSurfaceST\right)_{ij} \\ \quad \forall ij \in \EproST ,
\end{multline}

where conservation of energy on the network side $\Qproin = \Qprohx$ is defined in the model equations. In principle, part of the produced heat could be disposed or stored in a storage. We will, however, assume that all heat is supplied instantaneously to the DHN, as will be explained in section \ref{sec:stateConstraints}.

\subsection{Additional technological state constraints}
\label{sec:stateConstraints}
In addition to satisfying the physical model (flow and heat transfer equations), technological constraints $\defStateConstraints \leq 0$ and $\defStateConstraintsEquality = 0$ are defined to ensure that a useful optimization problem is solved. In this study, four sets of technological constraints are defined: heat demand satisfaction, $\STName$ heat incorporation, max pressure, and producer capacity constraints.
\subsubsection{Equality constraints}
\subsubsection*{Heat demand satisfaction}
The heat demand of each consumer in each period must be satisfied, this can be formulated as 
\begin{multline}\label{eq:ThermalComfort}
	\frac{\Qdemandi{\gi\gj,\timeVar}-\heat_{\gi\gj,\timeVar}}{\Qdemandi{\gi\gj,\timeVar}} = 0 \quad \forall ij \in \Econ ,
\end{multline}
where $\Qdemandi{\gi\gj,\timeVar}$ is the heat demand and $\heat_{\gi\gj,\timeVar}$ the transferred heat.

\subsubsection*{Solar thermal heat incorporation}
As mentioned in section \ref{section:STModel} we enforce full heat integration of the available solar heat into the DHN. This can be done via
\begin{multline}\label{eq:STHeat}
 \frac{\Qprosol - \Qprohx}{\maxCapacityIJ} = 0 \quad \forall ij \in \EproST ,
\end{multline}
with their respective definitions given in equations \ref{eq:HXST3} and \ref{eq:HXST2}.

\subsubsection{Inequality constraints}
\subsubsection*{Max Pressure}
A maximal pressure increase of $\maxPressure = \maxPressureValue \, \mathrm{bar}$ is defined on each producer edge which can be stated as
\begin{multline}\label{eq:MaxPressure}
	 \frac{(\pressure_{j,\timeVar}-\pressure_{i,\timeVar})-\maxPressure}{\maxPressure} \, \leq 0 \quad \forall ij \in \Epro .
\end{multline}

\subsubsection*{Producer capacity}
Finally, we limit the heat input of each producer in each period by its thermal capacity as
\begin{equation}\label{eq:capacityConstraint}
	\gijEdge{\left(\frac{\flow_\timeVar \Delta \dTinf_\timeVar}{\maxCapacity}\right)} \spHeatCap \rho - \gijEdge{\capVar} \leq 0 \quad \forall ij \in \Epro \setminus \EproST.
\end{equation}
Note that this constraint is not defined for the $\STName$ producers since we already ensure for those producers via equation \ref{eq:STHeat} the full incorporation of the available heat.

\subsection{A DHN model and optimization framework}
\label{DHNmodelandOpt}

\begin{figure}[h!]
	\centering
	\includegraphics[width=1.0\columnwidth]{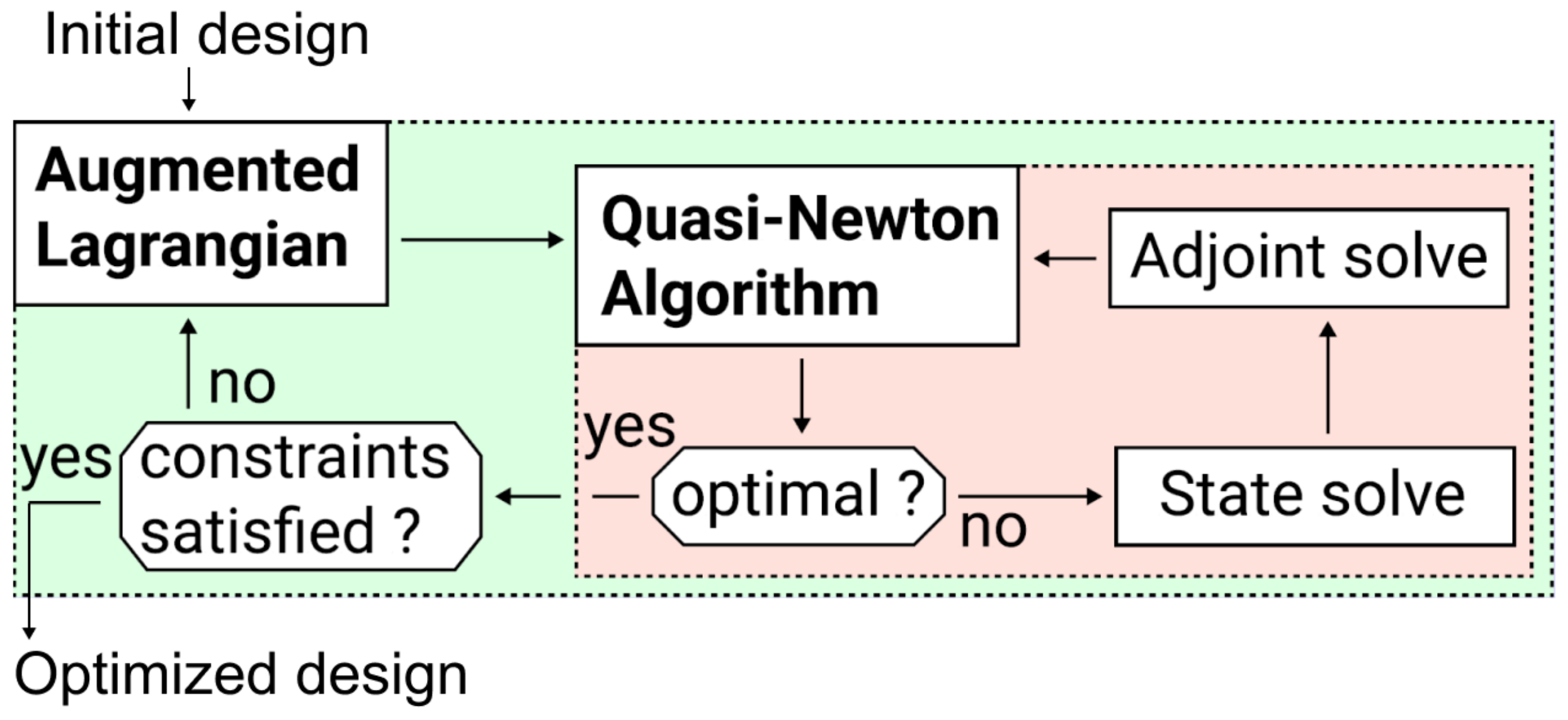}
	\caption{A flowchart of the optimization methodology. It consists of two nested components: The Augmented Lagrangian approach enforcing the technological (state) constraints and the Quasi-Newton algorithm solving the resulting bound-constrained non-linear sub-problem. Figure adapted from \citet{WackTopology}.}
	\label{fig:OptimizationRoutine}
\end{figure}

To minimize the costs and $\COTWO$ emissions of an existing DHN by optimizing the design and operation of its newly built producers, the physics of future network operation must be accurately captured to account for the distributions of temperature, pressure, and flow rate. Therefore, the optimization problem is constrained by a set of non-linear model equations  $\equalCon_{\timeVar}(\designVar_{\timeVar},\stateVar_{\timeVar})$ representing the hydraulic and thermal transport problem in the network. This network model ensures conservation of mass, momentum, and energy for the heat suppliers, consumers, and pipes in the network. 
The temporal scope of this problem is modeled in a multi-period approach assuming a steady state in each individual period\footnote{For a detailed discussion of this approach, the reader is referred to \citet{WackMP}}. Hence, neither hydraulic or thermal transients are considered. The complete network model and its equations can be found in \citet{WackMP}. 
\\
The non-linear optimization problem is solved using an Augmented Lagrangian approach combined with a Quasi-Newton method based on adjoint gradients, see figure \ref{fig:OptimizationRoutine}. Similar to the forward problem, the adjoint problem of each period can be solved independently, allowing parallel computation of the optimization problem across time periods \cite{WackMP}.

\subsection{Temporal resolution and time aggregation}
\begin{figure}[h!]
	\centering
	\includegraphics[width=\columnwidth]{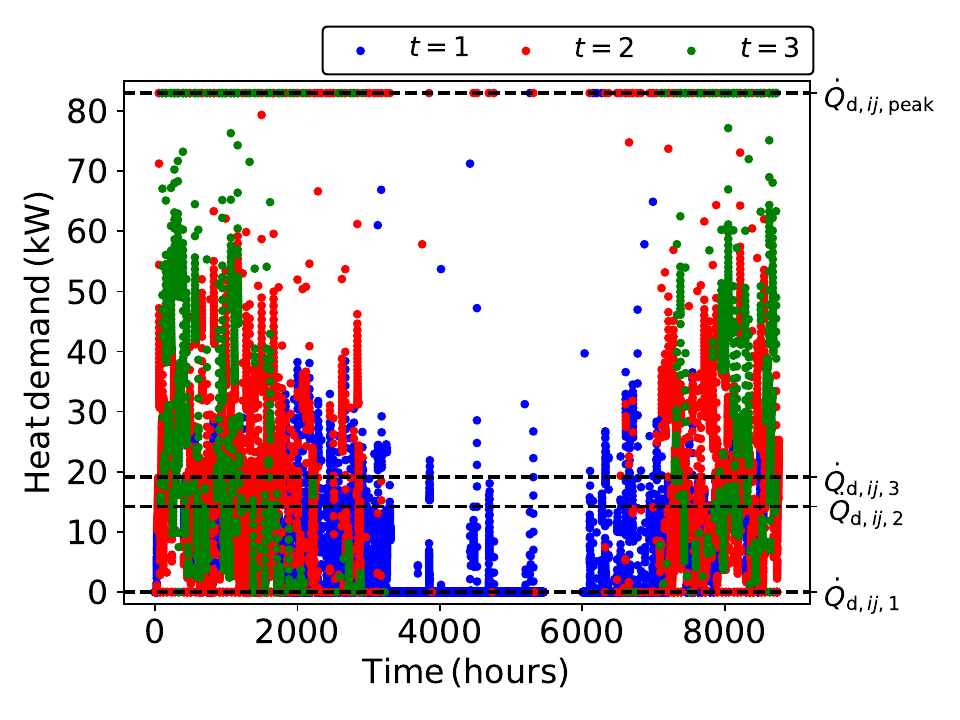}
	\caption{Exemplary heat demand $\Qdemandi{\gi\gj}$ of one of the considered consumers throughout the year and its aggregation into 3 representative periods plus 1 peak period. Figure adapted from \citet{WackMP}.}
	\label{fig:heat_demand}
\end{figure}

\begin{figure}[h!]
	\centering
	\includegraphics[width=\columnwidth]{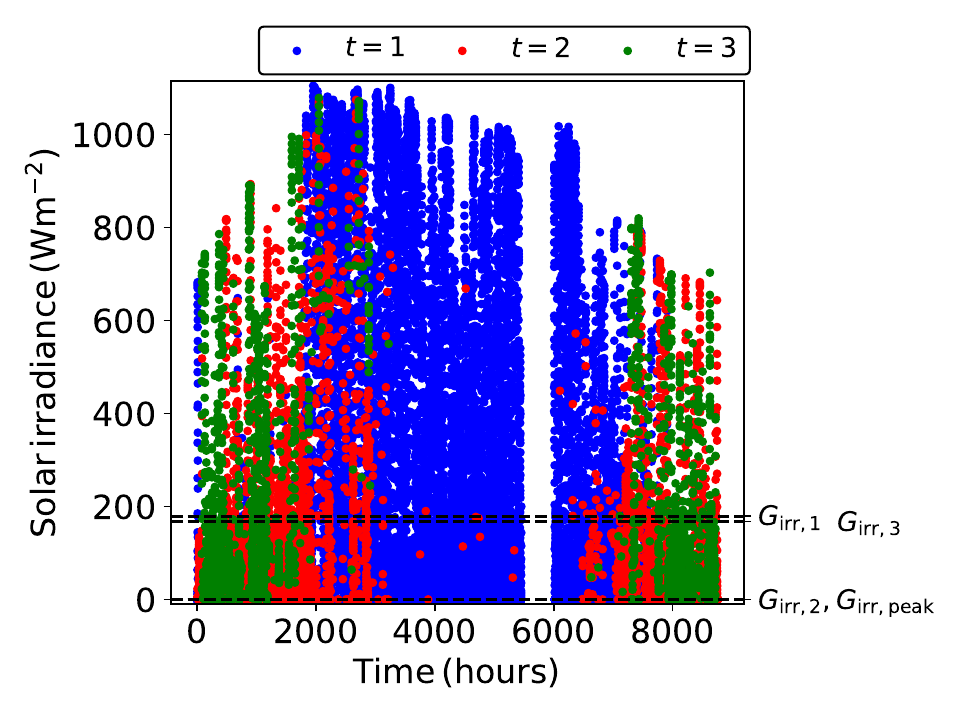}
	\caption{The solar irradiance $\Girr$ throughout the year and its aggregation into 3 representative periods plus 1 peak period.}
	\label{fig:irradiance}
\end{figure}

\begin{figure}[h!]
	\centering
	\includegraphics[width=\columnwidth]{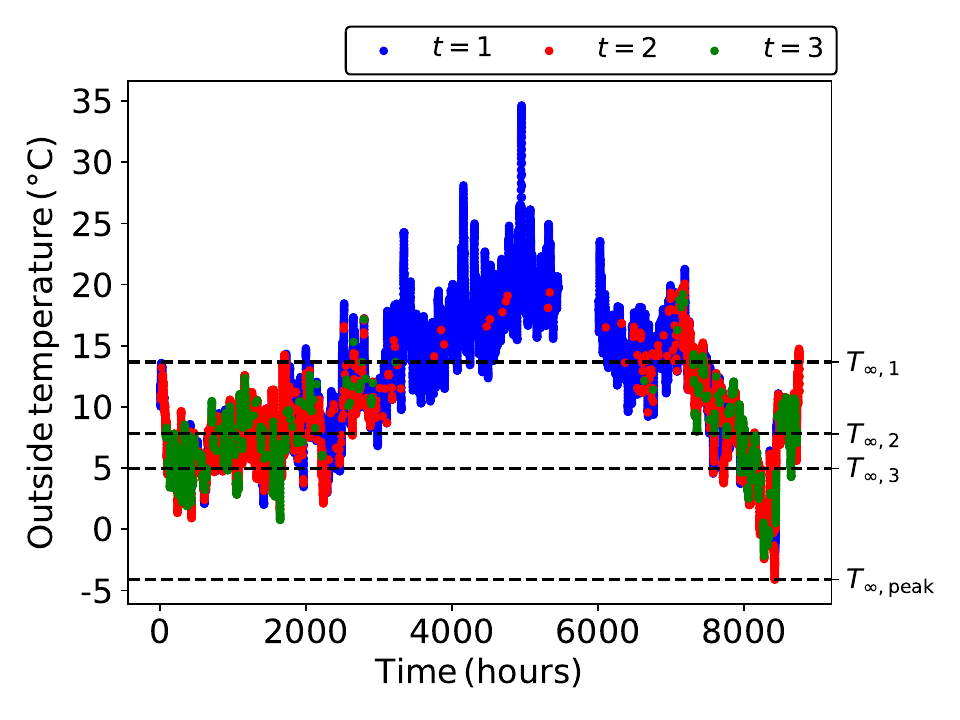}
	\caption{The outside temperature $\TOutside$ throughout the year and its aggregation into 3 representative periods plus 1 peak period. Figure adapted from \citet{WackMP}.}
	\label{fig:temperature}
\end{figure}

The multi-period approach in this study allows to consider a number of time-dependent parameters that significantly influence the producer design and operation of DHNs. In this work, the focus is on temporal variations of the consumer heat demands, the outside temperature and the solar irradiance. To control the computational complexity, the time-dependent parameters are temporally aggregated into a number of representative periods. In this work, we follow the same clustering approach as it was done in a previous work by the authors of this paper \cite{WackMP}, which is a  k-medoids clustering approach described by Kotzur et al. \cite{Kotzur2018}. Additionally, a consecutive series of time steps with no heat demands in the entire network (summer period) is excluded from the time series. The weight of the summer period of the solved case study in this work is 0.063 (6.3 \% of the year) which results in a conversion factor $\KOPEX$ (number of active operating hours per year) for the OPEX and $\COTWO$ costs of $\KOPEXvalue$ $\unit{\hour\per\year}$. Figures \ref{fig:heat_demand}, \ref{fig:irradiance}, and \ref{fig:temperature} provide a visualization of the time series data and the resulting aggregated periods.

\newcommand{\nAttributes}{\n_{\mathrm{at}}}
\newcommand{\nObservations}{\n_{\mathrm{ob}}}
\newcommand{\nCluster}{\n_{\mathrm{cl}}}
\newcommand{\veQdemand}[1][]{\ve{\heat}_{\mathrm{d}#1}}

To ensure feasibility of the optimized producer design under worst-case conditions, an additional worst-case period $\peakPeriod\in\setPeriods$ is added to the set of aggregated periods. In this period, the worst-case of each attribute, e.g. the peak demand of each consumer $\max_{\timeVar}(\Qdemandi{\gi\gj,\timeVar})$, the lowest outside temperature $\min_{\timeVar}\left(\TOutsideTime\right)$, and the lowest solar irradiance $\min_{\timeVar}\left(\GirrTime\right)$ are combined and added as an additional period to the optimization. Here, only the \gls{capex} cost and feasibility of this peak period is considered by setting the period weight to zero ($\periodWeight_{\peakPeriod}=0$).

\FloatBarrier
\section{A heat producer retrofit multi-objective optimization for heating networks - a case study}\label{chap:simple}

\begin{figure*}[h!]
	\centering
	\includegraphics[width=2.0\columnwidth]{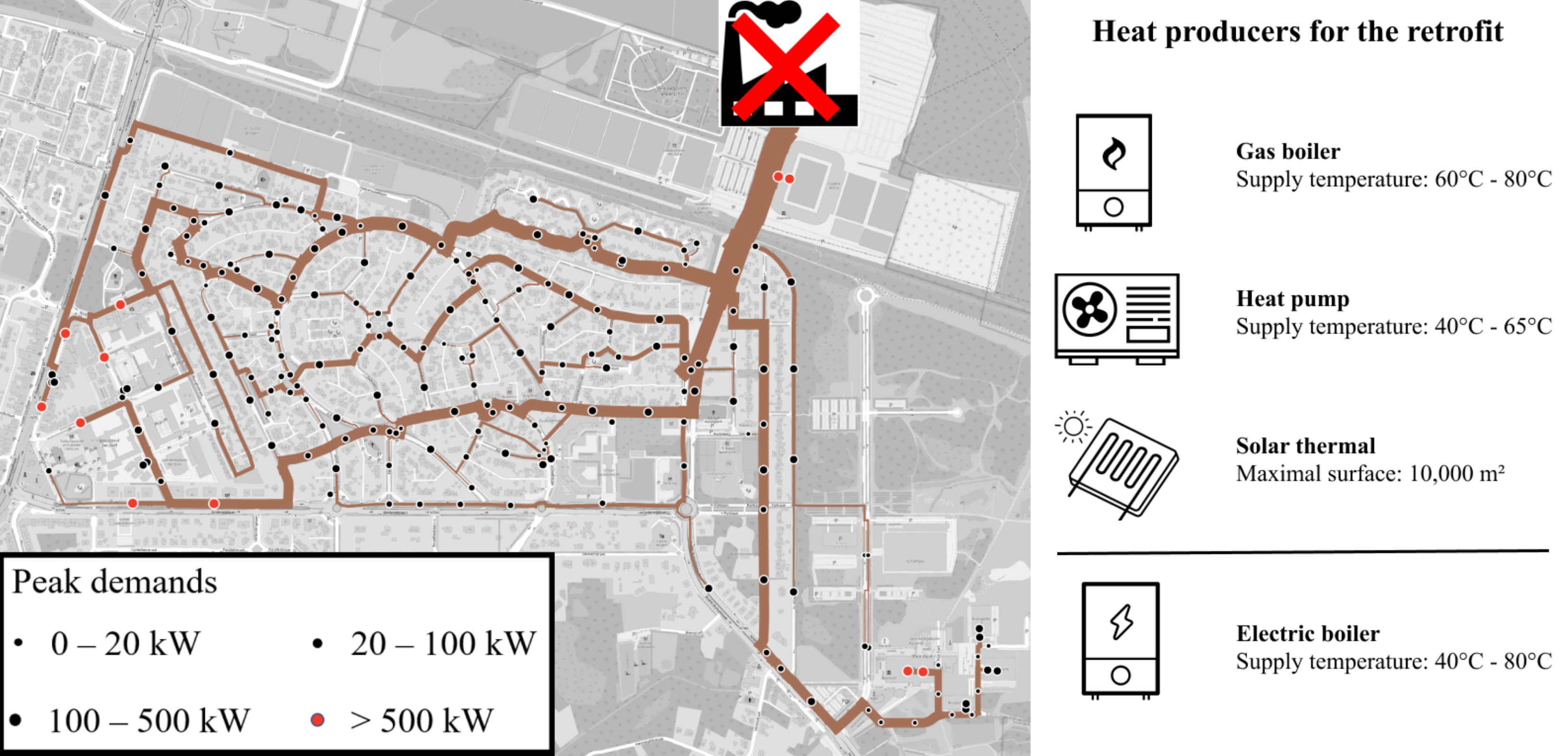}
	\caption{On the left: View of the Waterschei neighborhood in Genk, Belgium, highlighting the location of the assumed phased-out heat producer and the distribution of the 217 aggregated demand points. On the right: The available heat producers for the producer retrofit.}
	\label{fig:casesetup}
\end{figure*}

In this section, we assess the potential of an automated producer retrofit optimization approach, balancing network costs with $\COTWO$ emissions, based on multiple periods, for the producer retrofit of existing heating networks. We set up a fictitious case of an existing 3rd generation DHN in the Waterschei neighborhood of Genk, Belgium. The old heat source for this network must be phased out and replaced by new heat producers. Three heat producer types are available for the retrofit at first: A $\GBName$, a $\HPName$ and a $\STName$ unit; The impact of considering an $\EBName$ is investigated at a later stage as well. A visualization of the case study setup is provided in figure \ref{fig:casesetup}. Moreover, we assess the impact of different $\COTWO$ prices, ranging from 0 to 0.3\,$\si{\sieuro}\mathrm{\,kg^{-1}}$, on the producer design and operation. The neighborhood consists of 3800, mostly residential, buildings. The network topology and the pipe diameters were derived with a previous version of the design optimization framework developed by the authors of this paper \cite{WackTopology}. For this, a single high-temperature heat producer at 80°C was assumed in the north and only the peak heat demands were considered, optimizing for a worst-case design.

\definecolor{redNice}{HTML}{D88B84}
\definecolor{blueNice}{HTML}{8B84D8}
\definecolor{greenNice}{HTML}{84D88B}

\subsection*{Case setup}

The studied case is similiar to the \gls{dhn} development project in Waterschei used in \cite{WackMP}. Real heat demands are used as an input for the optimization. The locations of the heat demands were derived from \gls{gis} data (as described in \cite{Salenbien2023}). The total annual heat demand of the houses in the neighborhood is sourced from the database \textit{Warmtekaart Vlaanderen} \cite{Warmtekaart} and spatially aggregated per street segment into 217 representative demands. Annual variations of the heat demands were considered via representative annual heat demand curves based on the IDEAS library \cite{Jorissen2018} and the StROBe library \cite{Baetens2016}. The heat demand curves were modeled for different building types and renovation status (see \citet{WackMP} for more detailed information on the heat demand curves). An example of these heat demand time series is visualized in figure \ref{fig:heat_demand}. For the time-dependent outside temperature $\TOutside$, an hourly time series of temperatures for the year 2022 was obtained from the Royal Meteorological Institute of Belgium \cite{Temperature} for a weather station in Flanders. This time series is shown in figure \ref{fig:temperature}. For the time dependent solar irradiance $\Girr$ an hourly time series for the year 2020 was obtained from the radiation database PVGIS-SARAH2 \cite{IrradianceDataBase} for a city in Flanders (see figure \ref{fig:irradiance}). The slope (40°) and azimuth (-9°) of the surface were set in an optimal way by the PVGIS tool \cite{IrradianceTool} to maximize the usable irradiance.
For the $\GBName$, $\HPName$, and $\EBName$, minimal and maximal supply temperatures are defined which represent the lower and upper bound of the producer temperature design variable. Regarding their capacities, we assume that they can be built as large as needed. For the $\STName$ unit, the supply temperature is not optimized but is instead obtained by solving the heat exchanger equations between the $\STName$ loop and the DHN loop (see section \ref{section:STModel}). The maximally available surface area for the $\STName$ is bounded by the available area. The mentioned restrictions on the producer design can be found in figure \ref{fig:casesetup}.

After preprocessing and aggregation, the optimization problem has 217 heat consumers and $\timeSlices=4$ time periods, resulting in an optimization problem with around 891 design variables. The period weights $\periodWeight=\{0.65, 0.265, 0.085, 0\}$ together with the conversion factor $\KOPEX$ of $\KOPEXvalue$ $\si{\,\hour}\,\mathrm{yr^{-1}}$ amount to $\{222, 91, 29,0\}$ days of a representative year, respectively. The producer retrofit optimization problem for this case study is now solved using the proposed automated design approach. In the sections \ref{sec:CostsOnly} - \ref{sec:FixedSolar} only the $\GBName$, $\HPName$ and $\STName$ are considered. In section \ref{sec:ElectricBoiler} an $\EBName$ is added as well. 

\subsection{Economical design optimization}\label{sec:CostsOnly}
\begin{figure*}[h!]
	\centering
	\includegraphics[width=2.0\columnwidth]{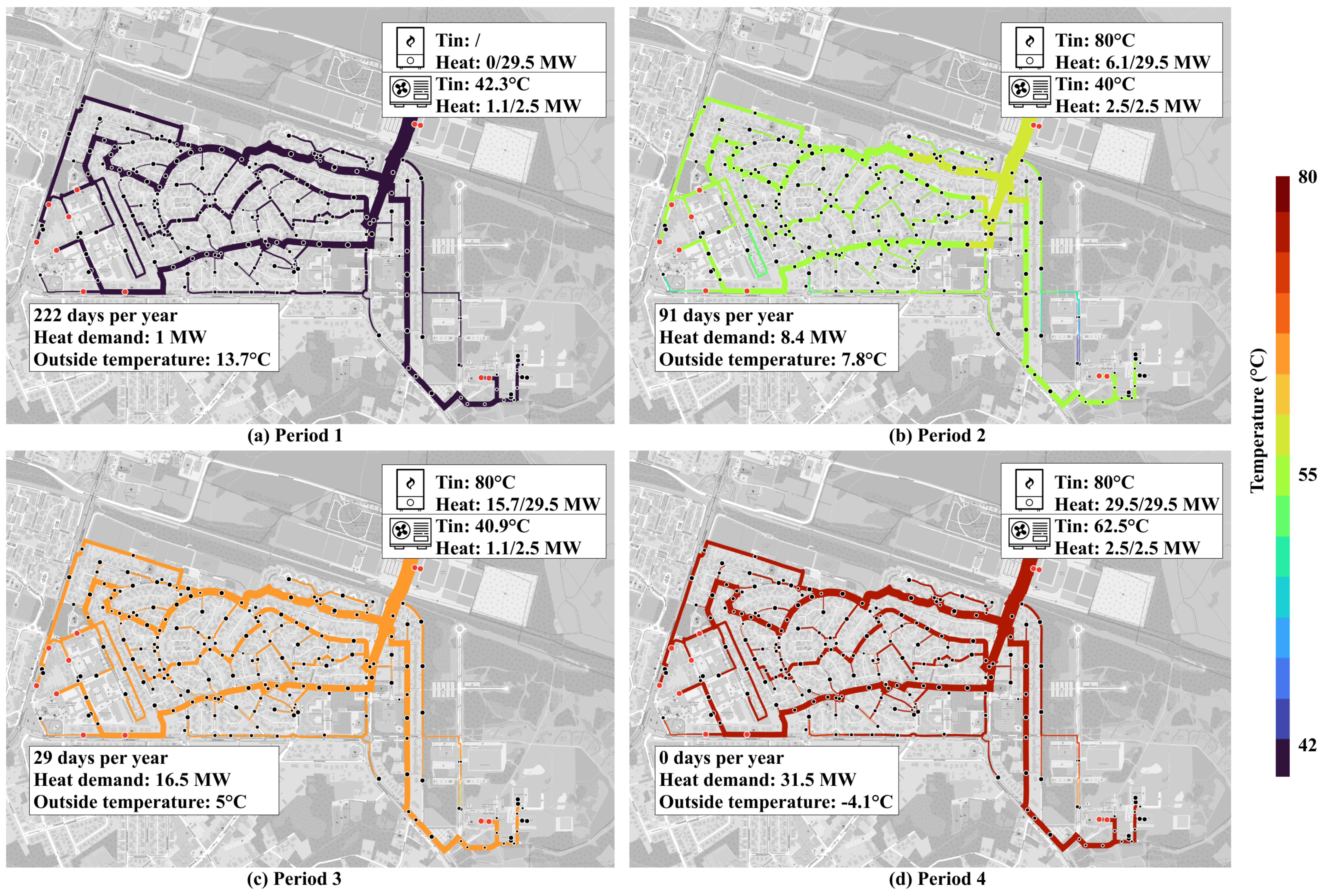}
	\caption{The result of the \textit{economical design optimization} showing the built producer capacities and their operation during each representative period. The network colors represent the network water temperatures.}
	\label{fig:CostsOnlyTemperature}
\end{figure*}

First, the heat producers for the DHN are optimized without the consideration of $\COTWO$ costs ($\COTWOPrice = 0$). The proposed optimization approach is used to determine for all producers the capacity $\capVar$ and inflow $\prodInput_{\timeVar}$, and for the $\GBName$ and $\HPName$ also the supply temperature $\prodTemp_{\timeVar}$. The result is summarized in figure \ref{fig:CostsOnlyTemperature}. The $\STName$ unit is not built at all due to a prohibitively low usable heat output and a high investment cost. This is further assessed in section \ref{sec:FixedSolar}. Moreover, the $\GBName$ is built with a capacity of 29.5$\si{\,\mega\watt}$ to cover the main heat load while the $\HPName$ is built with a capacity of only 2.5$\si{\,\mega\watt}$. The design choice is likely caused by the much higher specific investment cost of the $\HPName$ in comparison to the $\GBName$. Additionally, the lower efficiency of the $\GBName$ in comparison to the $\HPName$ is partially offset due to the low natural gas price which is cheaper by a factor 3 than the electricity. 

Regarding the operation, we see that 66\% of the total heat supply during the year (periods 1-3) comes from the $\GBName$, leading to a capacity factor of 9\%, while the operation of the $\HPName$ leads to capacity factor of 56\%, which is again likely caused by the different specific investment costs of the two units. 
The operational temperatures of the two units are driven by their respective efficiency curves. For the $\HPName$, a lower network supply temperature increases its COP and reduces thereby its electricity cost. Hence, we can see that the optimization lowers the $\HPName$ supply temperature towards its lower bound of $\si{40\,\degree\C}$ during the periods 1-3. In the peak period, however, its supply temperature is increased to 62.5$\si{\,\degree\C}$. This is on the one hand due to the pressure constraint that limits the flow rate. On the other hand, the temperature requirement of the consumers is with 55$\si{\,\degree\C}$ at its highest during the peak period, which requires a sufficiently high supply temperature, see \ref{app:ConsumerTemperatures} for more details on the temperature requirements of the consumers. For the $\GBName$, there is no price associated to its supply temperature. Hence, the highest possible temperature is selected in all periods because it allows for lower flow rates, thereby lowering pumping costs. 

In general, we can see that the overall supply temperature is increasing with increasing heat demand over the different periods to satisfy the heat demand, while meeting the limit on the maximal flow rate, imposed by the pressure constraint. Moreover, in period one, the longest period of the year, the $\GBName$ is completely shut down and all the heat is supplied by the operationally cost-efficient $\HPName$. For the periods 2-4 on the other hand, the two producer units are combined and their supplied hot water streams are mixed. The optimization thereby selects a hybrid heat pump operation as economically optimal.

\begin{figure}[h!]
	\centering
	\includegraphics[width=\columnwidth]{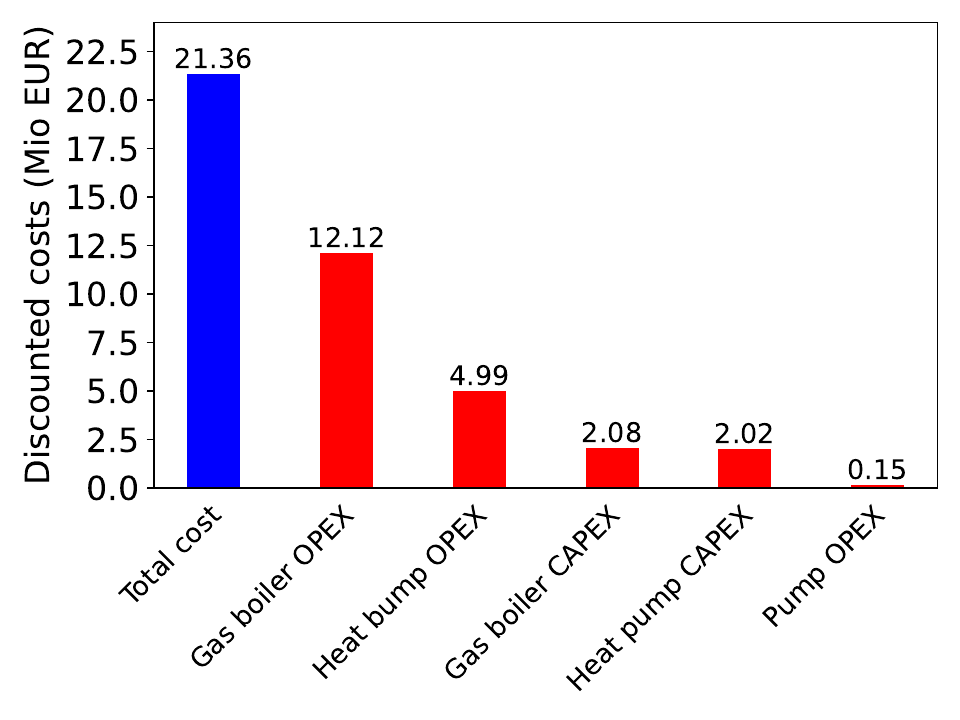}
	\caption{Cost analysis of the \textit{economical design optimization}.}
	\label{fig:caseII_costs}
\end{figure}

The analysis of the discounted costs in figure \ref{fig:caseII_costs} shows that the operational cost account for more than 75\% of the total cost of which more than two third are coming from the operation of the $\GBName$. The investment (CAPEX) costs of the heat producers are basically the same for the $\GBName$ and the $\HPName$ with around $\si{2\,\mega\sieuro}$ each, while the combined CAPEX amount to around 19\% of the total cost. The pumping cost represent less than 1\% of the total cost. 
	
\FloatBarrier
\subsection{Impact of $\COTWOITALIC$ pricing}
\label{sec:COTWOpricing}
\begin{figure}[h]
	
	\subfloat[Capacity]{
		\includegraphics[clip,width=\columnwidth]{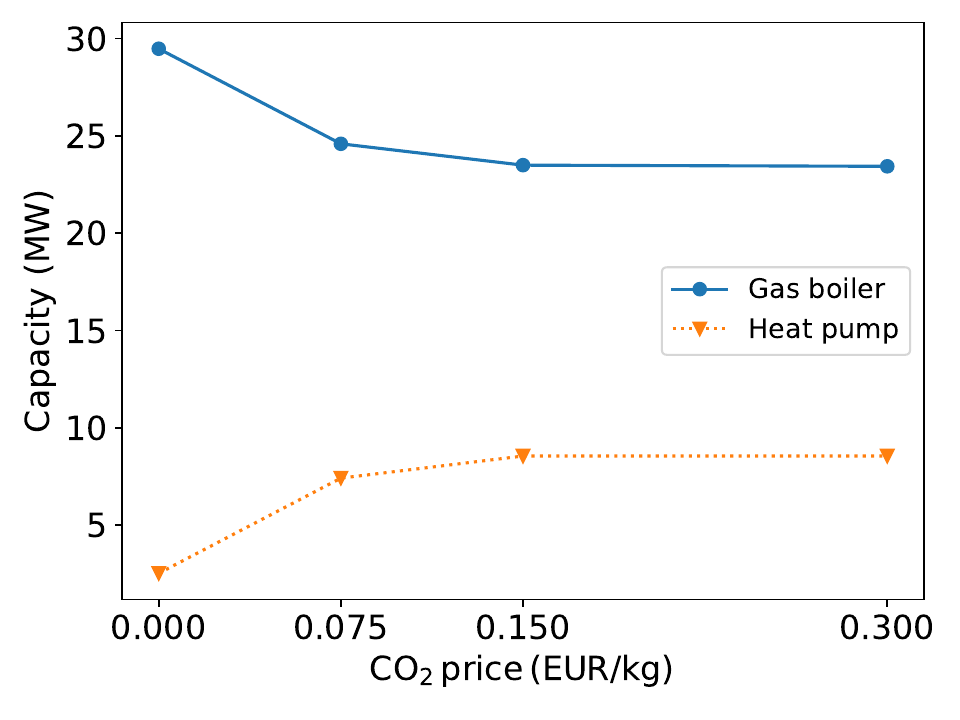}%
		\label{fig:CO2_capacities}
	}
	
	\subfloat[Heat supply share]{
		\includegraphics[clip,width=\columnwidth]{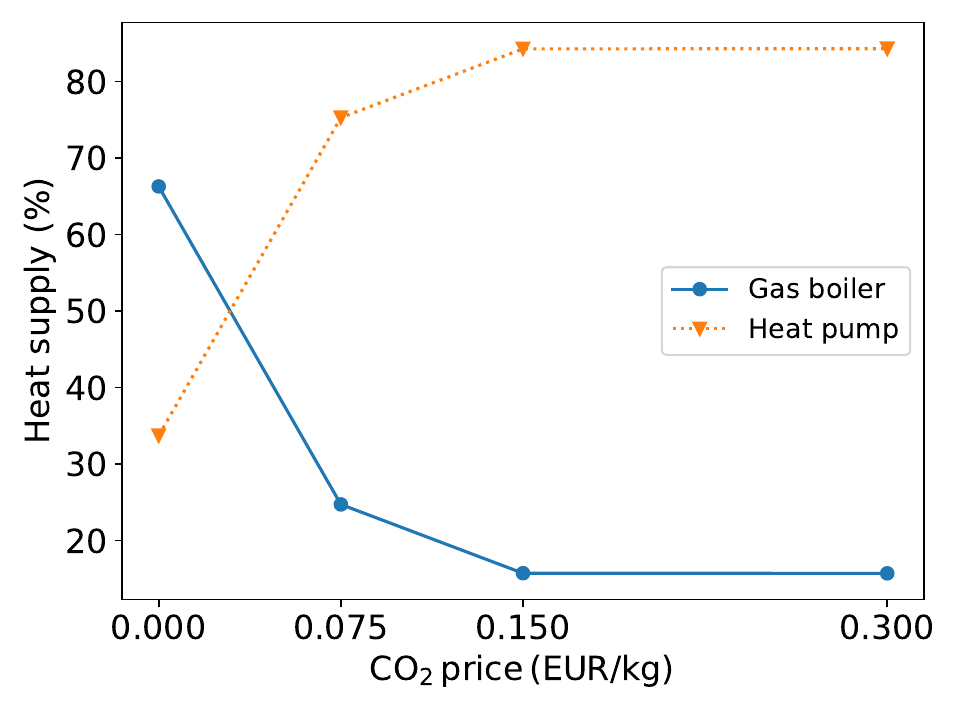}%
		\label{fig:CO2_HeatShare}
	}
	
	\caption{Changes in the (a) capacities and (b) shares of supplied heat of the heat producers $\HPName$ and $\GBName$ when changing the $\COTWO$ price. The result is from the \textit{impact of $CO_2$ pricing} case.}
	\label{fig:CO2_capacity_and_HeatShare}
\end{figure}

Now, the producer capacities, inflows, and for the $\GBName$ and $\HPName$ also the supply temperatures, are optimized while considering $\COTWO$ costs. To assess the impact of the $\COTWO$ price on the producer retrofit we consider three different $\COTWO$ prices in \,$\si{\sieuro}\mathrm{\,kg^{-1}}$: $\COTWOPrice \in \{0.075,0.15,0.3\}$. This range of $\COTWO$ prices is based on the Stern review \cite{Stern2007} which suggests that carbon pricing should be at a level that corresponds to the expected damage cost of using fossil fuels. 
\\
Analyzing the result, we see again that the $\STName$ unit is not built at all due to its unfavorable heat output - cost relation. Looking at the built capacities of $\HPName$ and $\GBName$ for the different $\COTWO$ prices in figure \ref{fig:CO2_capacities}, we can see that the $\HPName$ capacity increases with an increasing $\COTWO$ price from 2.5$\si{\,\mega\watt}$ up to 8.6$\si{\,\mega\watt}$. The opposite is true for the $\GBName$ where the capacity decreases from 29.5$\si{\,\mega\watt}$ to 23.4$\si{\,\mega\watt}$. If we analyze the share of heat supply from both units, we can see in figure \ref{fig:CO2_HeatShare} that the supply share from the $\HPName$ rises from 34\% up to 84\% while the share from the $\GBName$ decreases equivalently from 66\% to 16\%. Hence, despite the $\GBName$ capacity being still almost a factor 3 larger than the $\HPName$ capacity, the main heat load is supplied from the $\HPName$ in the case of the highest $\COTWO$ price. The same holds for the two other scenarios considering $\COTWO$ prices. The $\GBName$ is therefore operated more and more as a peak unit.

\begin{figure}[h]
	
	\subfloat[Discounted cost]{
		\includegraphics[clip,width=\columnwidth]{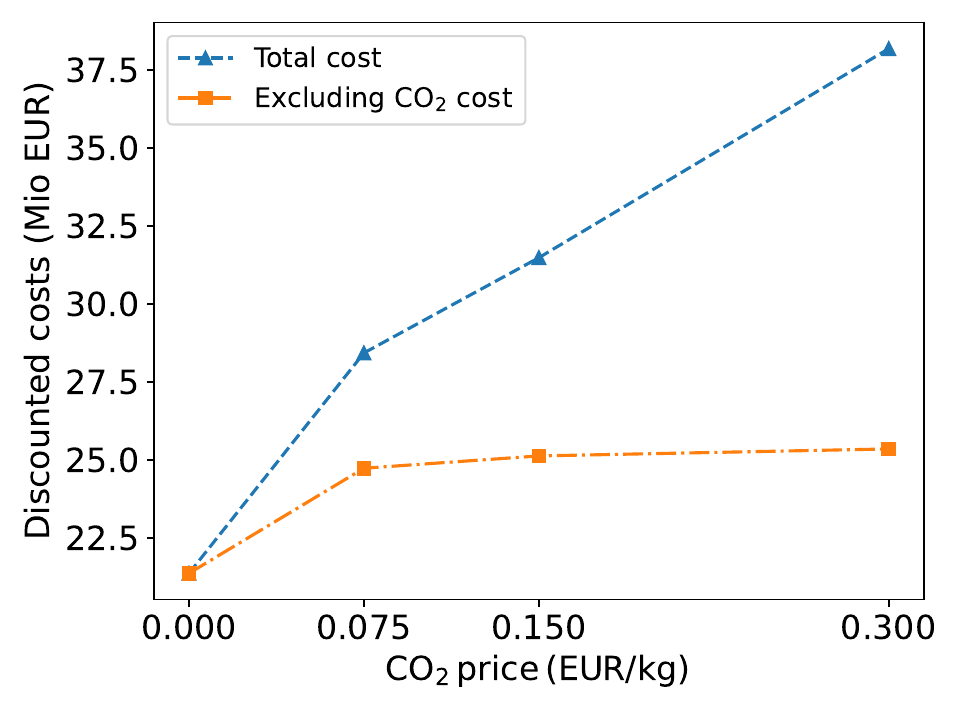}%
		\label{fig:CO2_NPVs}
	}
	
	\subfloat[Specific $\COTWO$ emissions per producer]{
		\includegraphics[clip,width=\columnwidth]{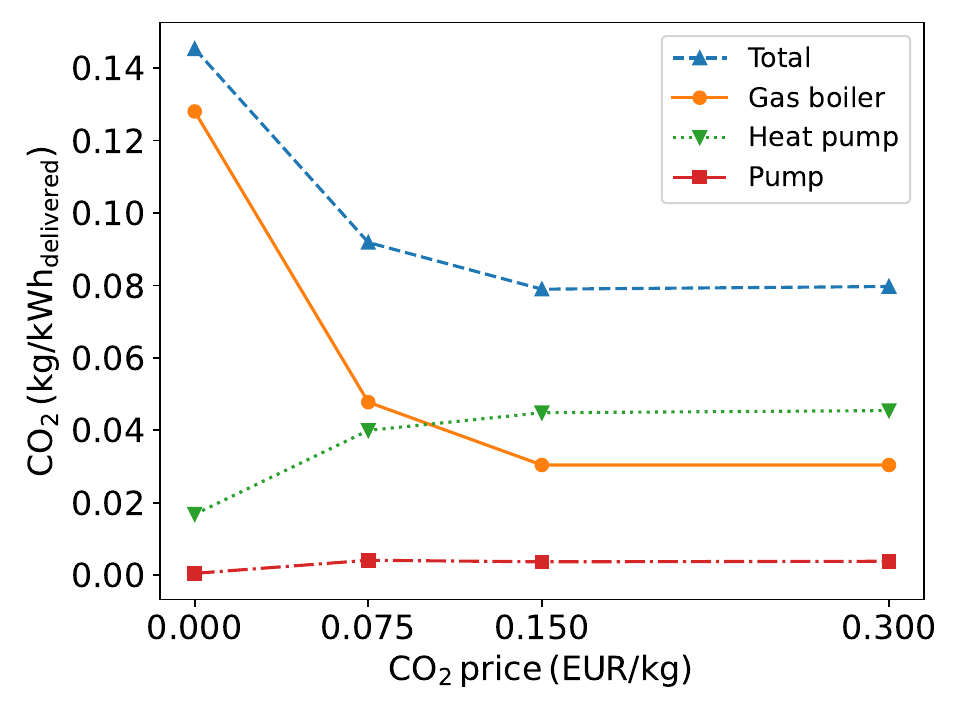}%
		\label{fig:CO2_SpecificCO2emissions_perProducer}
	}
	
	\caption{Changes in the (a) discounted cost and (b) specific $\COTWO$ emissions per producer of the $\HPName$ and $\GBName$ when changing the $\COTWO$ price. The result is from the \textit{impact of $CO_2$ pricing} case.}
	\label{fig:CO2_NPV_and_CO2}
\end{figure}

Moreover, it can be clearly seen in figures \ref{fig:CO2_capacities} and \ref{fig:CO2_HeatShare} that for the considered $\COTWO$ prices the changes in both the design (capacities) and operation (heat shares) stagnate and that the $\GBName$ is still built. To put the $\COTWO$ prices in perspective, the highest considered $\COTWO$ price of 0.3$\si{\,\sieuro}\,\mathrm{kg^{-1}}$ is almost by a factor of 3 higher than the highest $\COTWO$ price from the European Union Emission Trading System (EU-ETS) \cite{EUETS} in the year 2023 which was at 0.105$\si{\,\sieuro}\,\mathrm{kg^{-1}}$. 
The quasi-stagnation of the design and operation can be seen equivalently in the discounted project cost in figure \ref{fig:CO2_NPVs}. From a $\COTWO$ price of 0.075$\si{\,\sieuro}\,\mathrm{kg^{-1}}$ onward, the further price increase is primarily increasing the $\COTWO$ cost of the project. The total project cost increases by 78.8\% from 21.4$\si{\,\mega\sieuro}$ to 38.2$\si{\,\mega\sieuro}$. Regarding the $\COTWO$ emissions, we can see in figure \ref{fig:CO2_SpecificCO2emissions_perProducer} that the total $\COTWO$ emissions per $\si{\kilo\watt\hour}$ of delivered heat of the DHN project decrease from 0.145$\si{\,\kilogram\per{\kilo\watt\hour}}$ to 0.08$\si{\,\kilogram\per{\kilo\watt\hour}}$, a reduction of 45\%. Looking at the different sources of the $\COTWO$ emission, we can see that for the case of no $\COTWO$ price, the origin of the $\COTWO$ emissions is clearly dominated by the $\GBName$. For the $\COTWO$ prices of 0.15 and 0.3$\si{\,\sieuro}\,\mathrm{kg^{-1}}$, more $\COTWO$ is emitted by the $\HPName$ than by the $\GBName$. However, this must put in relation to their heat supply which is for the $\HPName$ for these two prices more than a factor of 4 larger than for the $\GBName$, corresponding to significantly lower specific $\COTWO$ emissions for the $\HPName$. Figure \ref{fig:CO2_SpecificCO2emissions_perProducer} also reveals that the $\COTWO$ contribution from the circulation pumps is almost negligible in comparison to the heat producers. For the $\COTWO$ price of 0.3$\si{\,\sieuro}\,\mathrm{kg^{-1}}$, for example, the $\COTWO$ emissions of the circulation pumps are a factor of 10 smaller than the $\COTWO$ emissions of the $\GBName$.

\FloatBarrier
\subsection{Increasing renewable ambition for the network retrofit} 
\label{sec:FixedSolar}
\begin{figure*}[h!]
	\centering
	\includegraphics[width=2.0\columnwidth]{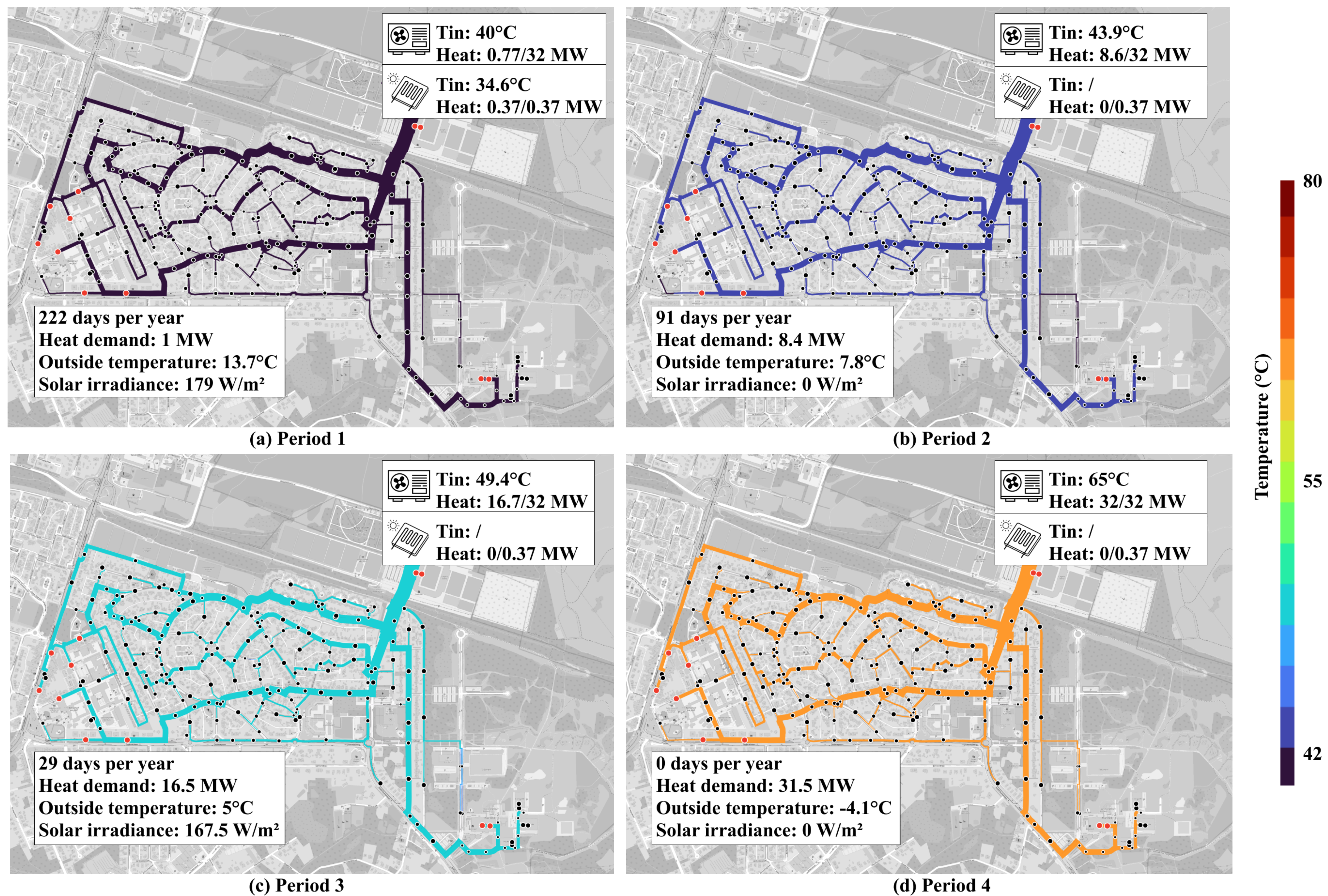}
	\caption{The result of the \textit{increasing renewable ambition for the network retrofit} case showing the producer capacities and their operation during each representative period. The network colors represent the network water temperatures.}
	\label{fig:caseIII_ST}
\end{figure*}

In the conducted case studies of sections \ref{sec:CostsOnly} and \ref{sec:COTWOpricing}, the $\STName$ unit is not build in the optimized heat producers design due to its prohibitive low usable heat output and high investment cost. To demonstrate nevertheless how a $\STName$ unit without a heat storage could contribute to the heat supply of a DHN and to assess its impact on the costs and $\COTWO$ emissions, we conduct a study in this section where we fix a priori the capacity of the $\STName$ unit to its maximal surface area of $10,000\si{\,\square\meter}$. Moreover, we exclude the $\GBName$ to focus on heat supply technologies which do not depend on fossil fuels and we set the $\COTWO$ price to 0. Hence, the remaining freedom of the optimization is given by the sizing of the $\HPName$ and the network operation. 
The result of this case study is summarized in figure \ref{fig:caseIII_ST}. Looking first at the solar irradiance of the aggregated periods, we can see that two of the periods (2 and 4) are defined with a solar irradiance of 0$\si{\,\watt\per{\square\meter}}$, making the operation of the $\STName$ during these periods impossible. Additionally, in the period 3, the combination of a low solar irradiance of 167.5$\si{\,\watt\per{\square\meter}}$ and a low outside temperature of 5$\si{\,\degree\C}$ leads to an average collector temperature $\TmeanST$ of around 31$\si{\,\degree\C}$ while some of the consumers have a return (cold) temperature of around 32$\si{\,\degree\C}$ in their heating system. This leads to a network return temperature that is too high to (fully) incorporate the solar heat. Hence, the $\STName$ unit is switched off in period 3 as well. 

Assessing the remaining period 1, where the $\STName$ unit is operated, we can see that it supplies around one third of the total heat supply or 0.37$\si{\,\mega\watt}$. Important to note is that the supply temperature on the network side of the $\STName$ unit is with a value of 34.6$\si{\,\degree\C}$ too low to satisfy consumer demands directly because they require a hot temperature of up to 36$\si{\,\degree\C}$ in their heating system in this period. Hence, even if the solar field would be larger, a temperature increase would still be needed as it is done here by mixing it with the hotter water flow from the $\HPName$.

\begin{figure}[h]
	\centering
	\includegraphics[width=1.0\columnwidth]{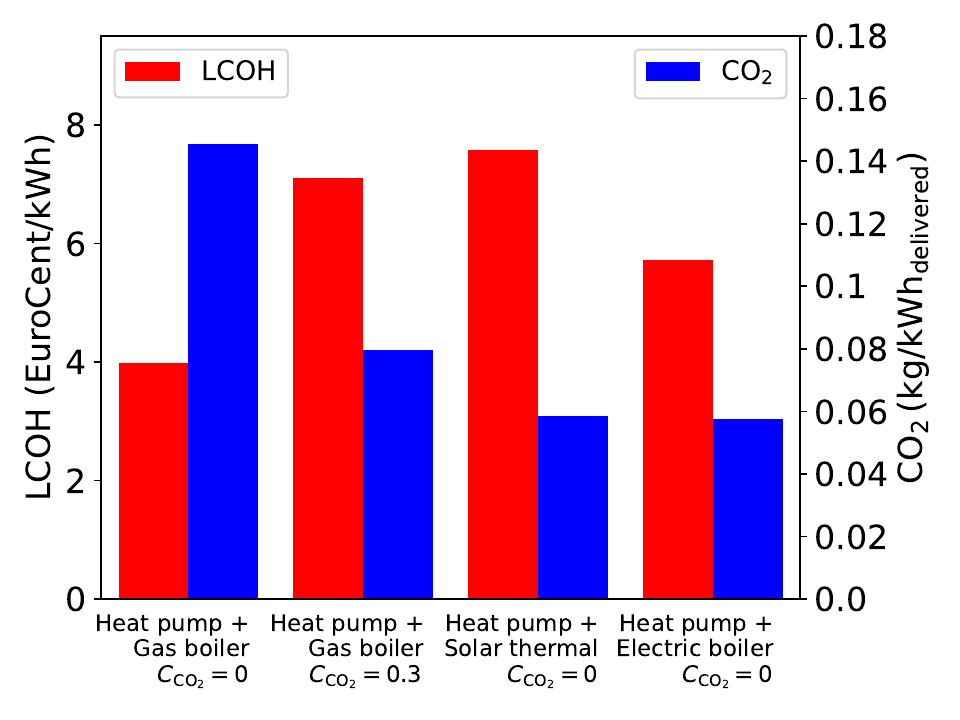}
	\caption{Levelized cost of heating (LCOH) and $\COTWO$ emissions per delivered $\si{\kilo\watt\hour}$ of heat for four of the optimized producer designs.} 
	\label{fig:four_cases_comparison_LCOH_CO2}
\end{figure}

We are now assessing the cost and $\COTWO$ emissions of this design and operation by comparing it with two of the  $\GBName$ + $\HPName$ designs from the previous section (the one without and the one with the highest $\COTWO$ price). We use for the cost comparison now the levelized cost of heating (LCOH), while the investment horizon of $\npvN = 30~\mathrm{years}$ and the discount rate of $\npvDiscount=5\%$ are unchanged. We can see in figure \ref{fig:four_cases_comparison_LCOH_CO2} that the obtained $\HPName$ + $\STName$ design is with a LCOH value of 0.076$\si{\,\sieuro\per{\kilo\watt\hour}}$ even more expensive than the $\GBName$ + $\HPName$ design with the highest considered $\COTWO$ price of 0.3$\si{\,\sieuro}\,\mathrm{kg^{-1}}$, explaining the $\STName$ phase out from the previous sections. These high cost are stemming from the necessity to build the (investment-wise) expensive $\HPName$ sufficiently large to provide heat to all consumers even under peak conditions. The specific $\COTWO$ emissions on the other hand are with a value of 0.058$\si{\,\kg\per{\kilo\watt\hour}}$ lower than the lowest $\COTWO$ emissions of a $\GBName$ + $\HPName$ combination at 0.08$\si{\,\kg\per{\kilo\watt\hour}}$. Both the cost and the $\COTWO$ emissions are mainly defined by the $\HPName$ and not by the $\STName$ unit since its heat contribution is of minor magnitude.

It becomes evident in this section that clustering the entire year into three periods is not sufficient to fully assess the potential of $\STName$ units for DHNs. Since the averaging of the irradiance and outside temperature essentially cuts out the most promising phases for solar heat production with high solar irradiance as it was highlighted in figures \ref{fig:ST_eta} and \ref{fig:ST_Tm} in section \ref{subsec:ST_model}. That the inclusion of intermittent renewables increases the temporal resolution required for a reliable economic assessment has also been highlighted by \citet{Merrick2016} and \citet{Pfenninger2017} for power system modeling.   
\FloatBarrier

\subsection{Electric boiler as an alternative peak boiler} 
\label{sec:ElectricBoiler}
\begin{figure*}[!htb]
	\centering
	\includegraphics[width=2.0\columnwidth]{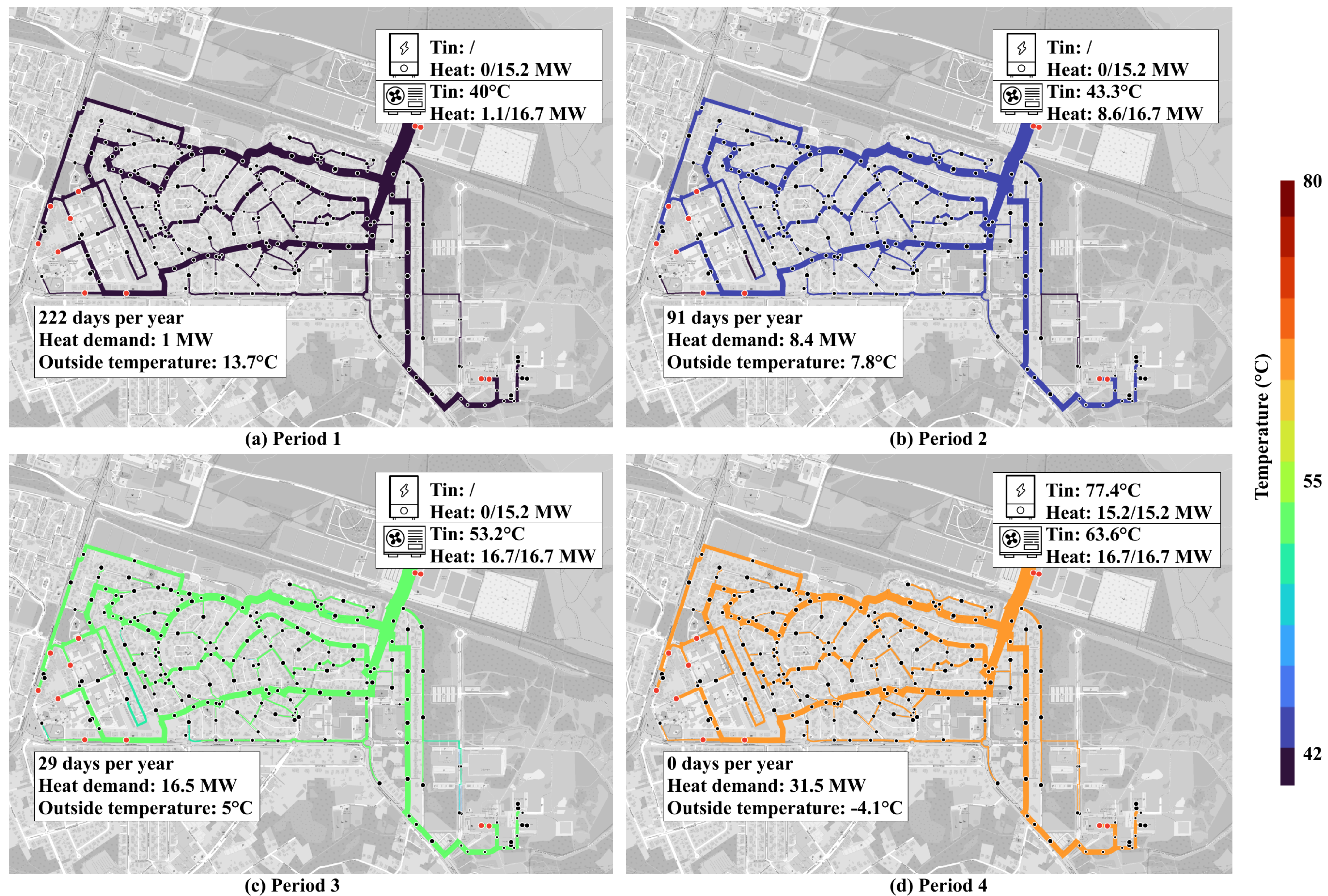}
	\caption{The result of the \textit{$\EBNameITALIC$ as an alternative peak boiler} case showing the producer capacities and their operation during each representative period. The network colors represent the network water temperatures.}
	\label{fig:caseII_EB_Temperature}
\end{figure*}

Now, we want to asses how an $\EBName$ could contribute to the heat supply of a DHN as an alternative to a $\GBName$. Electric boiler units have similar advantages over heat pumps as their natural gas based counterparts: they come with lower upfront investment cost and faster heat supply in particular for larger temperature differences between network return and supply. Hence, they are well suited to function as peak and back-up supply units without depending on fossil fuels. In this section, the producer capacities, inflows and the supply temperatures of the $\EBName$ and $\HPName$ are optimized while the $\COTWO$ price is set to zero ($\COTWOPrice = 0$). The $\GBName$ and the $\STName$ unit are excluded a priori, for the latter because sections \ref{sec:CostsOnly} and \ref{sec:COTWOpricing} showed that it is not competitive under the given circumstances. \\
The result of this case study is summarized in figure \ref{fig:caseII_EB_Temperature}. We can see that the $\HPName$ is built big enough to supply all the heat demand in the representative periods 1-3 (where operational costs occur) and that the $\EBName$ is not used at all in these three periods. The additional heat demand in the peak period is then covered by the $\EBName$. The difference in comparison to the $\GBName$ + $\HPName$ design is stemming from the higher investment cost of the $\EBName$ but even more from the usage of expensive electricity (versus cheap natural gas for the $\GBName$). Hence, the high efficiency of the $\HPName$ is exploited to counter the high costs of the only available fuel (electricity). An $\EBName$ of 15.2$\si{\,\mega\watt}$ is foreseen solely as peak unit. Looking at the cost of this design and operation in figure \ref{fig:four_cases_comparison_LCOH_CO2}, it can be seen that the $\EBName$ + $\HPName$ design is with LCOH of 0.057$\si{\,\sieuro\per{\kilo\watt\hour}}$ between the $\GBName$ + $\HPName$ design, without $\COTWO$ cost, and the $\HPName$ + $\STName$ design. Considering the $\COTWO$ emissions, it is positioned at the same level as the $\HPName$ + $\STName$ design with $\COTWO$ emissions of 0.058$\si{\,\kg\per{\kilo\watt\hour}}$. Hence, the $\EBName$ + $\HPName$ setup can be seen as a cost-effective producer design to decarbonize the heat supply.

\FloatBarrier
\section{Conclusions and outlook}
\subsection{Conclusions}
This paper proposes a mathematical optimization approach for the producer retrofit and decarbonization of existing DHNs based on a multi-period, physics-based, and multi-objective optimization problem, balancing CAPEX and OPEX with $\COTWO$ emission cost. The optimizer chooses the producer types, their capacities, and their supplied heat and supply temperatures. 

A case study for a medium-sized 3rd generation DHN is conducted and it is observed that the automated design approach was able to optimally retrofit the heat producers. For the study, a $\GBName$, a $\HPName$, and a $\STName$ unit are considered as potential heat producers. Four different $\COTWO$ prices are considered, ranging from 0 up to 0.3$\si{\,\sieuro}\,\mathrm{kg^{-1}}$. The results show that the $\COTWO$ pricing strongly impacts the design and the operation of the heat producers. The producer retrofit shifts toward a heat supply dominated by an energy-efficient and low-carbon $\HPName$ when the $\COTWO$ price is increased. The highest considered $\COTWO$ price resulted in a reduction of the specific $\COTWO$ emissions by 45\%, in comparison to the case with a $\COTWO$ price of 0. The LCOH on the other hand increased due to this $\COTWO$ price increase by 78.8\%. Moreover, it is observed that due to the high specific investment cost of a $\HPName$, the $\GBName$ is never completely phased out. Instead, it is still used to cover the high peak loads of the DHN due to its relatively low specific investment cost. Furthermore, the $\STName$ unit is never selected by the optimization as a part of the heat producer design due to its too low usable heat output in Belgium in the considered periods and its high investment cost.
Hence, from an economical point of view, a complete replacement of $\GBName$ units with $\HPName$ and $\STName$ units (without storage), by only using a $\COTWO$ pricing mechanism with prices of 0.3$\si{\,\sieuro}\,\mathrm{kg^{-1}}$ or less, is not achievable for the investigated case of a Belgian city. This result is affected by the current energy prices, the $\COTWO$ emission of the electricity mix, the heat producer investment costs, and by the selected periods as the resolution of the clustering impacts the usable heat supply from the $\STName$ producer.

Next, an $\EBName$ was assessed as an alternative to the $\GBName$. When no $\COTWO$ cost are considered, the setup $\HPName$ + $\EBName$ is able to reduce the $\COTWO$ emissions by 60\% to 0.058$\si{\,\kg\per{\kilo\watt\hour}}$ in comparison to a $\HPName$ + $\GBName$ setup. This $\COTWO$ reduction comes at an increase in LCOH of 44\% to 0.057$\si{\,\sieuro\per{\kilo\watt\hour}}$. Moreover, the LCOH of the $\HPName$ + $\EBName$ is 25\% cheaper in comparison to the $\HPName$ + $\STName$ setup while their $\COTWO$ emissions are at the same level. Hence, the $\HPName$ + $\EBName$ combination can be seen as one promising path to significantly lower the $\COTWO$ emissions of existing DHNs in a cost-efficient way. In all case setups it was clearly observed that the optimizer uses the freedom of setting the supply temperature to find the perfect balance between efficient heat generation, heat losses through the pipes, and pumping costs.
 
The automated, physics-based design tool presented in this work can be seen as an important prerequisite for the optimal producer retrofit of existing DHNs to transition from 2nd and 3rd to 4th generation DHNs. By balancing CAPEX and OPEX with $\COTWO$ emissions and defining the producer temperature as an additional degree of freedom, it provides a high design flexibility for DHN planners. It allows to analyze the impact of (potential) future $\COTWO$ prices on the network design, operation, and cost. Such an analysis can also be of high relevance for policy advisors and makers to assess the effectiveness of $\COTWO$ pricing mechanisms.
Moreover, the automated and scalable optimization approach allows to obtain results in short time and for large networks by only using a laptop\footnote{Each optimization ran less than 3 hours on a \emph{HP ZBook Power} mobile workstation with an \emph{intel i7} 14 core processor.}. The physics-based non-linear network simulation allows a reliable evaluation of the network operation with respect to the heat and pressure losses, the feasibility to satisfy heat demands and the obedience of pressure restrictions. Moreover, as introduced in this work, it allows to capture and eventually exploit temperature-dependent efficiencies of the heat producers.

\subsection{Outlook}

Future research should further investigate how a finer temporal resolution could help in better capturing phases of high solar irradiance and thereby higher yields of $\STName$ units. Moreover, including heat storage into the producer retrofit would enable a more efficient integration of $\STName$ units and would allow to reduce the needed peak capacity of the heat producers. In general, heat storage is expected to play a crucial role in maximizing renewable energy use in DHNs \cite{Jodeiri2022}. Moreover, it was observed that the peak period (the highest demand) strongly influences the design since it defines the needed heat producer capacity which leads to a preference of CAPEX-wise cheap (peak) heat producers and those are today classical fossil-fuel based producers, e.g., a $\GBName$. Hence, it is of interest to further investigate the peak demand requirements, by e.g. looking into simultaneity assessments of heat demands. Furthermore, as an alternative to the investigated $\COTWO$ pricing mechanism, strict $\COTWO$ emission constraints could be introduced to investigate the impact on design, operation, and cost of DHNs. 

Furthermore, a natural next research step would be to combine the presented heat producer optimization with the DHN topology optimization from \citet{WackTopology} for greenfield DHN projects. That this combination can be considered achievable from a computational point of view, is one of the main advantages of the proposed physics-based and scalable mathematical optimization methodology.

\section*{Data Availability}
A data-set including the structure, input parameters, time series and optimization results of the heating network and heat producers used in the case studies of this paper is available at the following link: \url{https://doi.org/10.48804/MNHGZX}. The optimization results can be replicated using the methodology and formulations described in this paper.

\section*{Acknowledgments}

Martin Sollich is funded by the Research Foundation – Flanders (FWO) through the PhD fellowship strategic basic research with the file number 1SH7624N.
Martin Sollich has received funding from the KU Leuven with the reference STG/21/016.

Yannick Wack has received funding from the Flemish institute for technological research (VITO).

\section*{CRediT authorship contribution statement}
\textbf{Martin Sollich}: Conceptualization, Data curation, Methodology, Software, Formal analysis, Investigation, Visualization, Writing – original draft, Funding acquisition. \textbf{Yannick Wack}: Conceptualization, Methodology, Software, Data curation, Writing – review \& editing. \textbf{Robbe Salenbien}: Software, Data Curation, Writing – review \& editing. \textbf{Maarten Blommaert}: Conceptualization, Methodology, Software, Supervision, Funding acquisition, Writing – review \& editing.

\section*{Compliance with ethical standards}
\subsubsection*{Conflict of interest}
The authors declare that they have no conflict of interest.
\subsubsection*{Funding}
The authors did not receive support from any organization for the submitted
work.

\appendix

\section{Additional details of the producer models}\label{app:model}

\subsection{Producer efficiencies}\label{app:efficiencies}

\begin{table}[h!]
	\centering
	\caption{The efficiency (curves) $\etaProdIJTime$ of the $\GBName$, $\HPName$, and $\EBName$ that link the supplied heat to the DHN with the energy input (natural gas or electricity). The parameter values are given in table \ref{tab:ProducerEfficienciesParameters}.}
	\label{tab:ProducerEfficiencies}
	\begin{tabularx}{\columnwidth}{>{\hsize=.65\hsize}X>{\hsize=.35\hsize}X}
		\toprule
		$\etaProdIJTimebold$ & \textbf{Producer} \\
		\midrule
		$\aGBeta\prodTempReturnCIJ^3+\bGBeta\prodTempReturnCIJ.^2+\cGBeta\prodTempReturnCIJ+\dGBeta$ & $\forall ij \in \EproGB$\\ 
		$\aHPeta*(\prodTempSupplyCIJ-\TOutsideTime)^2 + \bHPeta*(\prodTempSupplyCIJ-\TOutsideTime) + \cHPeta$ & $\forall ij \in \EproHP$\\
		$\aEBeta$ & $\forall ij \in \EproEB$\\
		\bottomrule
	\end{tabularx}
\end{table}

\begin{table}[h!]
	\centering
	\caption{The parameter values of the efficiency (curves) $\etaProdIJ$ of $\GBName$, $\HPName$, and $\EBName$ given in table \ref{tab:ProducerEfficiencies}.}
	\label{tab:ProducerEfficienciesParameters}
	\begin{tabularx}{\columnwidth}{>{\hsize=.4\hsize}X>{\hsize=.6\hsize}X}
		\toprule
		\textbf{Parameter} & \textbf{Value} \\
		\midrule
		$\aGBeta$ & $\aGBetavalue$ \\
		$\bGBeta$ & $\bGBetavalue$ \\
		$\cGBeta$ & $\cGBetavalue$ \\
		$\dGBeta$ & $\dGBetavalue$ \\ 
		$\aHPeta$ & $\aHPetavalue$ \\
		$\bHPeta$ & $\bHPetavalue$ \\
		$\cHPeta$ & $\cHPetavalue$ \\
		$\aEBeta$ & $\aEBetavalue$ \\
		\bottomrule
	\end{tabularx}
\end{table}

\begin{table}[h]
	\centering
	\caption{The parameter values of the solar collector efficiency curve given in equation \ref{eq:etaST}.}
	\label{tab:etaSTParameters}
	\begin{tabularx}{\columnwidth}{>{\hsize=.4\hsize}X>{\hsize=.6\hsize}X}
		\toprule
		\textbf{Parameter} & \textbf{Value} \\
		\midrule
		$\etaZeroST$ & $\etaZeroSTvalue$ \\
		$\aoneSTeta$ & $\aoneSTetavalue$ \\
		$\atwoSTeta$ & $\atwoSTetavalue$ \\
		\bottomrule
	\end{tabularx}
\end{table}

\FloatBarrier

\subsection{Producer CAPEX}\label{app:costs}
For each heat producer type, we fitted the specific investment cost based on literature and manufacturer data. The fits were obtained with the fit functionality in MATLAB \cite{fitMATLAB}.

\newcommand{\aGB}{a_{C,\GBIndex}}
\newcommand{\bGB}{b_{C,\GBIndex}}
\newcommand{\aHP}{a_{C,\HPIndex}}
\newcommand{\bHP}{b_{C,\HPIndex}}
\newcommand{\cHP}{c_{C,\HPIndex}}
\newcommand{\dHP}{d_{C,\HPIndex}}
\newcommand{\aEB}{a_{C,\EBIndex}}
\newcommand{\bEB}{b_{C,\EBIndex}}
\newcommand{\aST}{a_{C,\STIndex}}
\newcommand{\bST}{b_{C,\STIndex}}
\newcommand{\cST}{c_{C,\STIndex}}
\newcommand{\dST}{d_{C,\STIndex}}

\begin{table}[h!]
	\centering
	\caption{Fits of specific producer investment cost ($\prodCostSpecificIJ$). The producer capacity $\capVar*\maxCapacity$ is here indicated by $x$. The parameter values and references are provided in tables \ref{tab:SpecificInvestmentII} and \ref{tab:SpecificInvestmentIII}, respectively.}
	\label{tab:SpecificInvestmentI}
	\begin{tabularx}{\columnwidth}{>{\hsize=.5\hsize}X>{\hsize=.2\hsize}X>{\hsize=.3\hsize}X}
		\toprule
		$\prodCostSpecificIJbold$ & \textbf{Unit} & \textbf{Producer}\\
		\midrule
		$\aGB \, exp(\bGB \, x)$ & $\si{\sieuro}\mathrm{/W}$ & $\forall ij \in \EproGB$\\
		$\aHP \, exp(\bHP \, x)$ + $\cHP \, exp(\dHP \, x)$ & $\si{\sieuro}\mathrm{/W}$ & $\forall ij \in \EproHP$\\
		$\aST \, x^3 + \bST \, x^2 + \cST \, x + \dST$ & $\si{\sieuro}\mathrm{/m^2}$ & $\forall ij \in \EproST$\\
		$\aEB \, exp(\bEB \, x)$ & $\si{\sieuro}\mathrm{/W}$ & $\forall ij \in \EproEB$\\
		\bottomrule
	\end{tabularx}
\end{table}

\newcommand{\aGBValue}{0.1085}
\newcommand{\bGBValue}{-1.4578e-08}
\newcommand{\aHPValue}{0.3913}
\newcommand{\bHPValue}{-6.1444e-07}
\newcommand{\cHPValue}{0.7204}
\newcommand{\dHPValue}{-2.6479e-09}
\newcommand{\aEBValue}{0.19}
\newcommand{\bEBValue}{-2.0662e-08}
\newcommand{\aSTValue}{-4.8911e-13}
\newcommand{\bSTValue}{7.6116e-08}
\newcommand{\cSTValue}{-0.0041}
\newcommand{\dSTValue}{2.5209e+02}

\begin{table}[h!]
	\centering
	\caption{Parameter values of the specific investment cost $\prodCostSpecificIJ$ fits given in table \ref{tab:SpecificInvestmentI}.}
	\label{tab:SpecificInvestmentII}
	\begin{tabularx}{\columnwidth}{>{\hsize=.4\hsize}X>{\hsize=.6\hsize}X}
		\toprule
		\textbf{Fit parameter} & \textbf{Value} \\
		\midrule
		$\aGB$ & $\aGBValue$ \\ 
		$\bGB$ & $\bGBValue$ \\
		$\aHP$ & $\aHPValue$ \\
		$\bHP$ & $\bHPValue$ \\
		$\cHP$ & $\cHPValue$ \\
		$\dHP$ & $\dHPValue$ \\
		$\aST$ & $\aSTValue$ \\
		$\bST$ & $\bSTValue$ \\
		$\cST$ & $\cSTValue$ \\
		$\dST$ & $\dSTValue$ \\
		$\aEB$ & $\aEBValue$ \\
		$\bEB$ & $\bEBValue$ \\
		\bottomrule
	\end{tabularx}
\end{table}

\begin{table}[h!]
	\centering
	\caption{Sources of the specific investment costs of the producer models on which the fits given in tables \ref{tab:SpecificInvestmentI} and \ref{tab:SpecificInvestmentII} are based on.}
	\label{tab:SpecificInvestmentIII}
	\begin{tabularx}{\columnwidth}{XX}
		\toprule
		\textbf{Producer model} & \textbf{Source} \\
		\midrule
		Gas boiler & \cite{Heinen2016, Wirtz2020, Tommasi2018} \\
		Heat pump & \cite{Pieper2018} \\
		Solar thermal & \cite{Trier2018, LandPrice} \\
		Electric boiler & \cite{Akhtari2020, Hers2015} \\
		\bottomrule
	\end{tabularx}
\end{table}

\FloatBarrier

\subsection{Solar thermal model fits}

\begin{table*}[ht]
	\centering
	\caption{The fits of the hot $\temperature_{\secondaryIndex,j}$ and the cold $\temperature_{\secondaryIndex,i}$ $\STName$ temperature based on the measurement data of existing solar DHNs. The abbreviation "Temp." stands for Temperature. The parameter values are given in table \ref{tab:STtemperatureFitsParameters}.}
	\label{tab:STtemperatureFits}
	\begin{tabular}{lll}
		\toprule
		\textbf{Temp.} & \textbf{Fit} & \textbf{Producer} \\
		\midrule
		$\temperature_{\secondaryIndex,j,\timeVar} \, (\mathrm{hot})$ & $\aTjST + \bTjST\GirrTime + \cTjST\TOutsideTime + \dTjST\GirrTime^2 + \eTjST\GirrTime\TOutsideTime + \fTjST\TOutsideTime^2$ & $\forall ij \in \EproST$\\ 
		$\temperature_{\secondaryIndex,i,\timeVar} \, (\mathrm{cold})$ & $\aTiST + \bTiST\GirrTime + \cTiST\TOutsideTime + \dTiST\GirrTime^2 + \eTiST\GirrTime\TOutsideTime + \fTiST\TOutsideTime^2$ & $\forall ij \in \EproST$\\
		\bottomrule
	\end{tabular}
\end{table*}

\begin{table}[h]
	\centering
	\caption{The parameter values of the fits of the hot $\temperature_{\secondaryIndex,j}$ and the cold $\temperature_{\secondaryIndex,i}$ $\STName$ collector temperature given in table \ref{tab:STtemperatureFits}.}
	\label{tab:STtemperatureFitsParameters}
	\begin{tabularx}{\columnwidth}{>{\hsize=.4\hsize}X>{\hsize=.6\hsize}X}
		\toprule
		\textbf{Parameter} & \textbf{Value} \\
		\midrule
		$\aTjST$ & $\aTjSTvalue$ \\
		$\bTjST$ & $\bTjSTvalue$ \\
		$\cTjST$ & $\cTjSTvalue$ \\
		$\dTjST$ & $\dTjSTvalue$ \\ 
		$\eTjST$ & $\eTjSTvalue$ \\
		$\fTjST$ & $\fTjSTvalue$ \\
		$\aTiST$ & $\aTiSTvalue$ \\
		$\bTiST$ & $\bTiSTvalue$ \\
		$\cTiST$ & $\cTiSTvalue$ \\
		$\dTiST$ & $\dTiSTvalue$ \\ 
		$\eTiST$ & $\eTiSTvalue$ \\
		$\fTiST$ & $\fTiSTvalue$ \\
		\bottomrule
	\end{tabularx}
\end{table}

\FloatBarrier
\section{Consumer temperatures}\label{app:ConsumerTemperatures}

\begin{figure}[h]
	\centering
	\includegraphics[width=1.0\columnwidth]{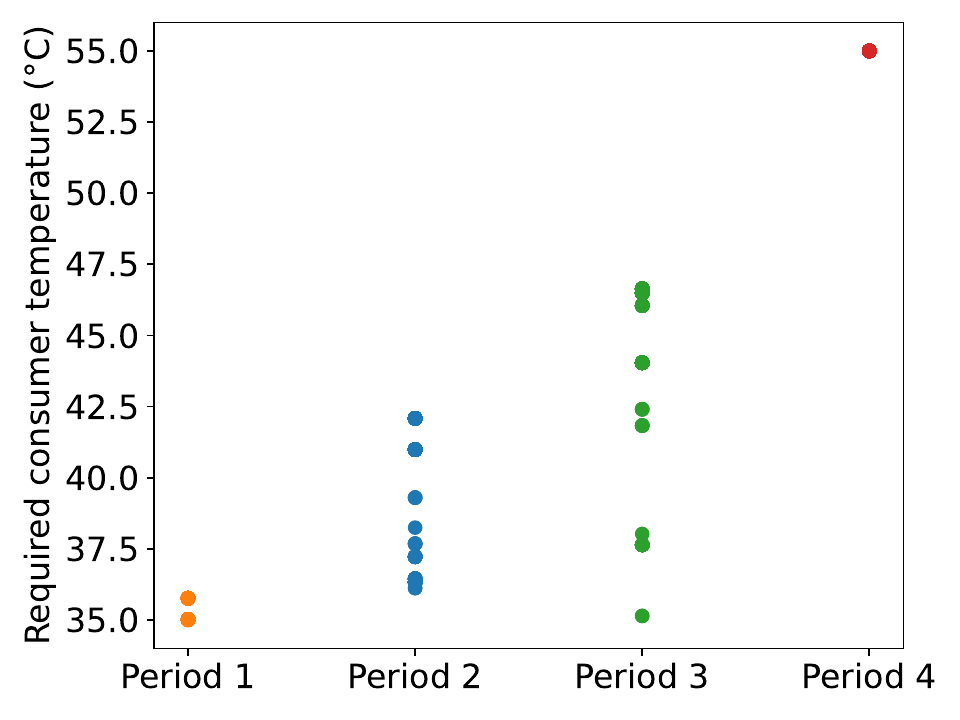}
	\caption{Required temperatures in the consumer heating systems (radiators) for each of the representative periods.}
	\label{fig:TconsHot}
\end{figure}

The full substation and consumer model can be found in \citet{WackMP}. Here, we show in figure \ref{fig:TconsHot} the required temperatures of the consumers within their heating system (secondary side) for each period. Note that we set the design temperature under peak conditions (highest demand) to 55$\si{\,\degree\C}$ for all consumers. The required temperature for the other periods was derived for each consumer individually based on its demand in that period and by using the characteristic radiator equation \cite{Feldhusen2014}. By modeling the dependency of the required consumer temperature on the heat demand, we allow the optimization to vary the DHN supply temperature throughout the year, enabling a more energy efficient heat supply and facilitating the integration of low-temperature renewable and waste heat sources. 

\FloatBarrier 

\bibliographystyle{elsarticle-num-names}
\bibliography{main}

\end{document}